\documentclass{mcom-l}

\usepackage{amssymb}
\usepackage{amsthm}
\usepackage{amsmath}
\usepackage{amsxtra}
\usepackage{mathrsfs}
\usepackage[all]{xy}

\theoremstyle{plain}
\newtheorem{thm}{Theorem}[section]

\newtheorem{prop}[thm]{Proposition}

\newtheorem{lem}[thm]{Lemma}
\newtheorem{cor}[thm]{Corollary}
\newtheorem*{thm*}{Theorem}
\newtheorem*{prop*}{Proposition}
\newtheorem*{cor*}{Corollary}

\theoremstyle{remark}

\newtheorem{exm}[thm]{Example}
\newtheorem{defn}[thm]{Definition}
\newtheorem{rmk}[thm]{Remark}
\newtheorem{notn}[thm]{Notation}

\newenvironment{enumalph}
{\begin{enumerate}}
{\end{enumerate}}

\newenvironment{enumroman}
{\begin{enumerate}}
{\end{enumerate}}

\DeclareMathOperator{\Cl}{Cl}
\DeclareMathOperator{\disc}{disc}
\DeclareMathOperator{\Frob}{Frob}
\DeclareMathOperator{\Gal}{Gal}
\DeclareMathOperator{\Hom}{Hom}

\DeclareMathOperator{\Pic}{Pic}
\DeclareMathOperator{\repart}{Re}
\DeclareMathOperator{\res}{res}

\newcommand{\bfC}{\mathbb C}

\newcommand{\bfQ}{\mathbb Q}

\newcommand{\bfZ}{\mathbb Z}

\newcommand{\fraka}{\mathfrak{a}}

\newcommand{\frakp}{\mathfrak{p}}

\newcommand{\calD}{\mathcal{D}}
\newcommand{\calP}{\mathcal{P}}

\newcommand{\la}{\langle}
\newcommand{\ra}{\rangle}
\newcommand{\eps}{\epsilon}

\newcommand{\Kbar}{\overline{K}}
\newcommand{\bfQbar}{\overline{\bfQ}}

\newcommand{\psmod}[1]{~(\textup{\text{mod}}~{#1})}
\newcommand{\psod}[1]{~({#1})}
\newcommand{\legen}[2]{\left(\frac{#1}{#2}\right)}

\newcommand{\dhat}{\widehat{d}}
\newcommand{\sbmin}{{}_{\textup{min}}}

\begin{document}

\title{Quadratic forms that represent almost the same primes}

\author{John Voight}
\address{Department of Mathematics, University of California, Berkeley, Berkeley, California 94720}
\email{jvoight@math.berkeley.edu}
\thanks{The author's research was partially supported by an NSF Graduate Fellowship.  The author would like to thank Hendrik Lenstra, Hans Parshall, Peter Stevenhagen, and the reviewer for their helpful comments, as well as William Stein and the MECCAH cluster for computer time.}


\date{September 16, 2005}


\begin{abstract}
Jagy and Kaplansky exhibited a table of $68$ pairs of positive definite binary quadratic forms that represent the same odd primes and conjectured that their list is complete outside of ``trivial'' pairs.  In this article, we confirm their conjecture, and in fact find all pairs of such forms that represent the same primes outside of a finite set.
\end{abstract}

\maketitle

\section{Introduction}

The forms $x^2+9y^2$ and $x^2+12y^2$ represent the same set of prime numbers, namely, those primes $p$ which can be written $p=12n+1$ for some positive integer $n$.  What other like pairs of forms exist?  Jagy and Kaplansky \cite{JagyKaplansky} performed a computer search for pairs that represent the same set of odd primes and found certain ``trivial'' pairs which occur infinitely often and listed other sporadic examples.  They conjecture that their list is complete.

Using the tools of class field theory, in this article we give a provably complete list of such pairs.  By a \emph{form} $Q$ we mean an integral positive definite binary quadratic form $Q=ax^2+bxy+cy^2 \in \bfZ[x,y]$; the \emph{discriminant} of $Q$ is $b^2-4ac=D=df^2<0$, where $d$ is the discriminant of $\bfQ(\sqrt{D})$ or the \emph{fundamental discriminant}, and $f \geq 1$.  We will often abbreviate $Q=\la a,b,c \ra$.  

Throughout, we look for forms that represent the same primes outside of a finite set---we say then that they represent \emph{almost the same primes}.  A form represents the same primes as any equivalent form under the action of the group $GL_2(\bfZ)$.  Hence from now on (except in the statement of Proposition \ref{quadforms}, see Remark \ref{SL2remark}, and in the proof of Lemma \ref{char2powerforms}), we insist that a form be \emph{$GL_2(\bfZ)$-reduced}, i.e., $0 \leq b \leq a \leq c$.  Moreover, the set of primes represented by a form is finite (up to a finite set, it is empty) if and only if the form is nonprimitive, that is to say $\gcd(a,b,c)>1$, and any two nonprimitive forms represent almost the same primes.  We therefore also insist that a form be \emph{primitive}, so that the set of primes represented is infinite.

If $Q_1,Q_2$ are forms which represent almost the same primes, we write $Q_1 \sim Q_2$; it is clear that $\sim$ defines an equivalence relation on the set of forms.  To every equivalence class $C$ of forms, we associate the set $\delta(C)$ of fundamental discriminants $d$ of the forms in $C$ as well as the set $\Delta(C)$ of discriminants $D$ of forms in $C$.  

The main result of this article is the following (Theorem \ref{thm1}).

\begin{thm*}
There are exactly\/ $67$ equivalence classes $C$ of forms with $\# \delta(C) \geq 2$.  There are exactly\/ $6$ classes with $\#\delta(C)=3$ and there is no class with $\# \delta(C) \geq 4$.  
\end{thm*}

\begin{cor*}
There are exactly\/ $111$ pairs of forms $Q_1,Q_2$ with fundamental discriminants $d_1 \neq d_2$such that $Q_1 \sim Q_2$.
\end{cor*}

The forms are listed in Tables \textup{1--5} at the end of this article.

As a complement to this theorem, we characterize forms $Q_1 \sim Q_2$ with the same fundamental discriminant $d_1=d_2$ (Theorem \ref{k1eqk2}).

\begin{thm*}
Let $Q_1=\la a_1,b_1,c_1 \ra$ be a form with $|D_1|>4$.  Then there exists a form $Q_2 \sim Q_1$ such that $|D_2|>|D_1|$ and $d_1=d_2=d$ if and only if one of the following holds:
\begin{enumroman}
\item $d \equiv 1 \pmod{8}$ and $2 \nmid D_1$;
\item $2 \mid D_1$ and either $b_1=a_1$ or $a_1=c_1$.
\end{enumroman}
\end{thm*}

These theorems together prove the conjecture of Jagy and Kaplansky in the affirmative regarding pairs that represent the same odd primes.  (See also Remark \ref{jkmiss} at the end of this article.)

We now give an outline of the proof.  To a form $Q$, we associate an ideal class in an imaginary quadratic order and, by the Artin map, to this ideal class we associate an element of a ring class group (Proposition \ref{quadforms}), and the representability of a prime $p$ by the form $Q$ then amounts to a certain splitting condition on $p$ in the ring class field associated to $Q$.  Therefore, two forms $Q_1,Q_2$ represent almost the same primes if and only if they give rise to the same splitting data, which can be formally thought of as an open and closed subset $S \subset \Gal(\bfQbar/\bfQ)$ (Lemma \ref{Pinj}).  By Galois theory, such a set has a (unique) minimal field of definition $L$ (Proposition \ref{minflddefexists}).  

We first treat the case when the forms $Q_1,Q_2$ have different fundamental discriminants $d_1 \neq d_2$.
Group theoretic considerations show that $Q_1,Q_2$ have the same genus class field, contained in the field $L$, and that their ring class groups are \emph{of type dividing $(2,\dots,2,4)$}, i.e.~they can be embedded in $(\bfZ/2\bfZ)^r \oplus \bfZ/4\bfZ$ for some $r \in \bfZ_{\geq 0}$ (Proposition \ref{thmclgrp}).  We then extend existing methods for bounding class groups of imaginary quadratic fields and, using a computer, effectively determine all possible ring class extensions which may arise from the forms $Q_1 \sim Q_2$ (\S 5).  From this finite data we can then list all possible pairs of quadratic forms which represent almost the same primes (\S 6).  

When $Q_1,Q_2$ have the same fundamental discriminant $d_1=d_2$, we can by classical methods determine necessary and sufficient conditions for $Q_1 \sim Q_2$ (\S 7).

As a side result which may be of independent interest, we provide the following classification of class groups of quadratic orders (Proposition \ref{alldiscrim}).  

\begin{prop*}
There are at least $226$ and at most $227$ fundamental discriminants $D=d<0$ such that $\Cl(d)$ is of type dividing $(2,\dots,2,4)$, and there are at least $199$ and at most $205$ such discriminants $D$ of nonmaximal orders.  
\end{prop*}

These orders are listed in Tables 7--16 at the end of this article.

\section{Ring class fields}

In this section, we fix notation and summarize without proof the few results we will need from class field theory and the theory of $L$-functions (see e.g. \cite{Cox}, \cite{Lang}, and \cite{Washington}).

Let $K=\bfQ(\sqrt{d})$ be an imaginary quadratic field of discriminant $d<0$ with ring of integers $A$. For an integer $f \geq 1$, consider the order $A_f=\bfZ+fA$; the discriminant of $A_f$ is $D=df^2$.  There is a bijection between the set $I(A)$ of ideals of $A$ coprime to $f$ and the set $I(A_f)$ of ideals of $A_f$ coprime to $f$, given by $\fraka \mapsto \fraka \cap A_f$ and conversely $\fraka_f  \mapsto \fraka_f A$.  Let $\Cl_f(d)=\Cl(D)=\Pic(A_f)$ be the class group of the order $A_f$, namely the group of invertible $A_f$-ideals modulo principal $A_f$-ideals.  Given an ideal $\fraka \subset A$ prime to $f$, the $A_f$-module $\fraka \cap A_f$ is trivial in $\Cl(D)$ if and only if $\fraka$ is principal and generated by an element $\alpha$ with $\alpha \equiv z \pmod{fA}$ for some $z \in \bfZ$.  We write $h_f(d)=h(D)=\# \Cl(D)$.  

\begin{prop}[{\cite[\S 9]{Cox}}] \label{ringclassfield}
There is a unique field $R_{(f)} \supset K$ inside $\overline{K}$ that is abelian over $K$ with the following properties:
\begin{enumroman}
\item Each prime $\frakp$ of $K$ coprime to $f$ is unramified in $R_{(f)}$;
\item There is an isomorphism
\begin{align*} 
\Cl_f(d) &\cong \Gal(R_{(f)}/K) \\
[\frakp \cap A_f] &\mapsto \Frob_\frakp
\end{align*}
\end{enumroman}
for each prime $\frakp$ of $K$ coprime to $f$.  

The field $R_{(f)}$ is the largest abelian extension of $K$ of conductor dividing $(f)$ in which all but finitely many primes of $K$ inert over $\bfQ$ split completely.  

The exact sequence
\[ 1 \to \Gal(R_{(f)}/K) \to \Gal(R_{(f)}/\bfQ) \to \Gal(K/\bfQ) \to 1 \]
splits, and a choice of splitting gives an isomorphism
\[ \Gal(R_{(f)}/\bfQ) \cong \Gal(R_{(f)}/K) \rtimes \Gal(K/\bfQ) \]
where the nontrivial element of $\Gal(K/\bfQ) \cong \bfZ/2\bfZ$ acts on $\Gal(R_{(f)}/K)$ by inversion $\sigma \mapsto \sigma^{-1}$.
\end{prop}

The field $R_{(f)}$ is called the \emph{ring class field of $K$ of modulus $f$}, and the map $\Cl_f(d) \cong \Gal(R_{(f)}/K)$ is known as the \emph{Artin isomorphism}.

\begin{rmk} \label{sigmasigmainverse}
As $\Gal(R_{(f)}/K)$ is abelian, we see from the proposition that the conjugacy class of an element $\sigma \in \Gal(R_{(f)}/K)$ in $\Gal(R_{(f)}/\bfQ)$ is equal to $\{\sigma,\sigma^{-1}\}$.  
\end{rmk}

\begin{cor} \label{rgcd}
Let $f_1,f_2 \in \bfZ_{\geq 1}$, and let $f=\gcd(f_1,f_2)$.  Then $R_{(f_1)} \cap R_{(f_2)} = R_{(f)}$. 
\end{cor}

\begin{proof}
The conductor of $R_{(f_1)} \cap R_{(f_2)}$ divides both $(f_1)$ and $(f_2)$, therefore it divides $(f)$ and has all but finitely many primes of $K$ inert over $\bfQ$ split completely, hence $R_{(f_1)} \cap R_{(f_2)} \subset R_{(f)}$.  Note also that $R_{(f)} \subset R_{(f_1)} \cap R_{(f_2)}$ since $f \mid f_1$ and $f \mid f_2$, therefore equality holds.
\end{proof}

\begin{prop}[{\cite[Theorem 7.7]{Cox}}] \label{quadforms}
Let $D=df^2<0$ be a discriminant.  Then there is a bijection between the set of $SL_2(\bfZ)$-reduced forms of discriminant $D$ and the set of ideal classes in $\Cl(D)$ by the identifications 
\[ Q=\la a,b,c \ra= ax^2+bxy+cy^2 \longleftrightarrow [\fraka]=[(a,(-b+f\sqrt{d})/2)] . \]

Let $Q$ be a form, with $Q \leftrightarrow [\fraka]$ for $\fraka$ an ideal of $A$ and $[\fraka]$ associated to $\sigma \in \Gal(R_{(f)}/K)$ under the Artin map.  Let $p \nmid f$ be prime.  Then $p$ is represented by $Q$ if and only if $[\fraka]$ contains an integral ideal of norm $p$, which holds if and only if we have $\Frob_p=\{\sigma,\sigma^{-1}\} \subset \Gal(R_{(f)}/\bfQ)$.  
\end{prop}

\begin{rmk} \label{SL2remark}
When considering primes represented by a form, we naturally link a Frobenius element together with its inverse; note that exactly one element of any conjugacy class $\{\sigma,\sigma^{-1}\}$ is associated with a $GL_2(\bfZ)$-reduced form.
\end{rmk}

\begin{rmk} \label{densityrmk}
Since $h_f(d)=[R_{(f)}:K]$, it follows from the Chebotarev density theorem that the density of the set of primes represented by $Q$ is equal to $1/(2h_f(d))$ if the corresponding element $\sigma$ has order $\leq 2$ (i.e., $\sigma=\sigma^{-1}$) and $1/h_f(d)$ otherwise.
\end{rmk}

\begin{lem} \label{repsame}
The forms $Q_1,Q_2$ represent almost the same primes ($Q_1 \sim Q_2$) if and only if for almost all primes $p$ of $\bfQ$, we have 
\[ \Frob_{p}=\{\sigma_1,\sigma_1^{-1}\} \subset \Gal(R_1/\bfQ) \Longleftrightarrow \Frob_{p}=\{\sigma_2,\sigma_2^{-1}\} \subset \Gal(R_2/\bfQ). \]
\end{lem}

\begin{rmk} \label{D1D2}
It follows from this that if $Q_1,Q_2$ are forms with the same discriminant $D_1=D_2$, then $Q_1 \sim Q_2$ if and only if $Q_1=Q_2$.
\end{rmk}

\begin{prop} \label{genusclassfield}
The field $P_{(f)} \subset R_{(f)}$ given by 
\[ \Gal(P_{(f)}/K) \cong \Cl_f(d)/\Cl_f(d)^2 \] 
is the largest subextension of $R_{(f)} \supset K$ with Galois group $\Gal(R_{(f)}/K)$ of exponent dividing $2$.  Moreover, the extension $P_{(f)} \supset \bfQ$ is itself abelian and of exponent $2$, and is the largest abelian subextension of $R_{(f)} \supset \bfQ$.
\end{prop}

The field $P_{(f)}$ is called the \emph{genus class field of $K$ of modulus $f$}.  

\begin{proof}
This follows immediately from Proposition \ref{ringclassfield}, as inversion acts trivially on a group of exponent dividing $2$.  
\end{proof}

We can compute the genus class field explicitly as follows.  For $p$ an odd prime we write $p^*=(-1)^{(p-1)/2}p$.  

\begin{cor} \label{computegenus}
Let $p_1,\dots,p_r$ be the odd primes dividing $D$ and let 
\[ K^*=K(\sqrt{p_1^*}, \dots, \sqrt{p_r^*}). \]
Then the genus class field $P_{(f)}$ of $K$ is as follows:
\[ P_{(f)}=
\begin{cases}
K^*(\sqrt{-1}), &\text{if }d \equiv 1 \psmod{4}\text{ and }4 \parallel f; \\
K^*(\sqrt{-1},\sqrt{2}), &\text{if }d \equiv 1 \psmod{4}\text{ and }8 \mid f; \\
K^*(\sqrt{2}), &\text{if }d \equiv 4 \psmod{8}\text{ and }4 \mid f; \\
K^*(\sqrt{-1}), &\text{if }d \equiv 0 \psmod{8}\text{ and }2 \mid f; \\
K^*, &\text{otherwise}.
\end{cases} \]
\end{cor}

\begin{proof}
See \cite[\S 6A]{Cox} for the case $f=1$.  The case $f>1$ is a standard calculation and follows in a similar way.
\end{proof}

\begin{cor} \label{samecor}
The odd primes $p$ which ramify in $P_{(f)}$ are exactly the odd primes that divide $D$.
\end{cor}

If $G$ is an abelian group and $n \in \bfZ_{>0}$, then we define $G[n]=\{g \in G:n g=0\}$.  

\begin{cor} \label{2gm1}
If $d$ has $g$ distinct prime factors, then $\Cl(d)[2] \cong (\bfZ/2\bfZ)^{g-1}$.
\end{cor}

For a fundamental discriminant $d<0$, let
\[ \chi(n)=\chi_{d}(n)=\left(\frac{d}{n}\right) \]
denote the Kronecker symbol. 

\begin{lem}[{\cite[Theorem 7.24]{Cox}}] \label{exactseq}
The sequence
\[ 1 \to A_f^* \to A^* \to (A/fA)^*/(\bfZ/f\bfZ)^* \to \Cl_f(d) \to \Cl(d) \to 1, \]
is exact, and
\[ h(D)=\frac{h(d)f}{[A^*:A_f^*]}\prod_{p \mid f}\left(1-\legen{d}{p}\frac{1}{p}\right). \]
\end{lem}

In the sequel, we will use lower bounds on the sizes of the class groups of quadratic fields.  If we write
\[ L(s,\chi)=\sum_{n=1}^{\infty} \frac{\chi(n)}{n^s}=\sum_{n=1}^{\infty}\frac{(d/n)}{n^s} \]
for $s \in \bfC$ with $\repart(s)>0$, then 
\[ h(d)=\frac{\sqrt{|d|}}{\pi}L(1,\chi) \]
for $|d|>4$ (see e.g. \cite[\S 6]{Davenport}).  By the Brauer-Siegel theorem, $\log h(d)$ is asymptotic to $\log(\sqrt{|d|})$ as $|d| \to \infty$; by a result of Siegel \cite{Siegel}, we know that for every $\eps>0$, there exists a constant $c(\eps)$ such that
\[ L(1,\chi) > \frac{c(\epsilon)}{|d|^\eps}; \]
however, this constant $c(\epsilon)$ is not known to be effectively computable.  Therefore we will use the following result on the size of $L(1,\chi)$.

\begin{lem}[Tatuzawa \cite{Tatuzawa}] \label{tatuzawa}
For any $0<\epsilon<1/2$, there is at most one fundamental discriminant $d<0$ with $\log |d|>\max(1/\epsilon,11.2)$ satisfying 
\[ L(1,\chi) \leq 0.655\frac{\epsilon}{|d|^{\eps}}. \]
\end{lem}

\section{Fields of definition}

We now proceed with a bit of Galois theory.  The reader may prefer on a first reading to skip to the next section and refer back when needed.  

Let $K$ be a field with separable closure $\Kbar$ and absolute Galois group $G=\Gal(\Kbar/K)$, equipped with the Krull topology.  Let $E$ be a finite extension of $K$ contained in $\Kbar$ and let $\Hom_K(E,\Kbar)$ denote the set of $K$-embeddings $E \hookrightarrow \Kbar$; if $E$ is Galois over $K$, then $\Hom_K(E,\Kbar)$ is identified with $\Gal(E/K)$.  We have a restriction map
\begin{align*} 
\res_E:G &\to \Hom_K(E,\Kbar) \\
\sigma &\mapsto \res_{E}(\sigma)=\sigma|_E.
\end{align*}
The map $\res_E$ is continuous if the finite set $\Hom_K(E,\Kbar)$ is equipped with the discrete topology.

\begin{lem} \label{flddefexists}
A subset $S \subset G$ is open and closed if and only if there exist a finite extension $L \supset K$ contained in $\Kbar$ and a set $T \subset \Hom_K(L,\Kbar)$ such that $S=\res_L^{-1}(T)$.
\end{lem}

\begin{proof}
Given $T \subset \Hom_K(L,\Kbar)$, note that $T$ is open and closed (in the discrete topology) and $\res_L$ is a continuous map.

Conversely, suppose $S \subset G$ is open and closed.  Then for every $\sigma \in S$, there exists an open neighborhood $U_\sigma=\res_{E_\sigma}^{-1}(\sigma|_{E_\sigma}) \subset S$ of $\sigma$ given by some finite extension $E_\sigma \supset K$.  Together these give an open cover $\{U_\sigma\}_{\sigma \in S}$ of $S$.  Since $G$ is compact and $S$ is closed, $S$ is itself compact and therefore is covered by $\{U_\sigma\}_{\sigma \in S'}$ for $S' \subset S$ a finite subset.  Let $L$ be the compositum of the fields $E_\sigma$ for $\sigma \in S'$, and let
\[ T=\{\tau \in \Hom_K(L,\Kbar):\tau|_{E_\sigma}=\sigma|_{E_\sigma}\text{ for some $\sigma\in S'$}\}. \]
Then by construction $S=\res_L^{-1}(T)$.
\end{proof}

\begin{defn}
Given an open and closed set $S \subset G$, we say that $L$ is a \emph{field of definition for $S$} if $L \supset K$ is a finite extension and there is a subset $T \subset \Hom_K(L,\Kbar)$ such that $S=\res_L^{-1}(T)$.
\end{defn}

\begin{rmk} \label{rmkflddef}
If $L$ is a field of definition with $S=\res_L^{-1}(T)$ for some subset $T \subset \Hom_K(L,\Kbar)$, then in fact $T=S|_L$.  Therefore $L$ is a field of definition for $S$ if and only if $\res_L^{-1}(S|_L)=S$, i.e.~for every $\sigma \in G$ and $\tau \in S$ such that $\sigma|_L=\tau|_L$ we have $\sigma \in S$.  It follows immediately from this that if $L$ is a field of definition for $S$ and $M \supset L$ is a finite extension, then $M$ is also a field of definition for $S$.
\end{rmk}

Put in these terms, Lemma \ref{flddefexists} states that every open and closed subset $S \subset G$ has a field of definition.  

\begin{defn}
A field of definition $L$ for $S$ is \emph{minimal} if for every field of definition $E$ for $S$, we have $L \subset E$.
\end{defn}

If a minimal field of definition $L$ exists, it is obviously unique.

\begin{prop} \label{minflddefexists}
For any open and closed set $S \subset G$, there exists a minimal field of definition $L(S)$ for $S$.  
\end{prop}

\begin{proof}
Consider the set
\[ H(S)=\{\sigma \in G:S\sigma=S\} \subset G; \] 
we claim that $L(S)=\Kbar^{H(S)}$.  

The set $H(S)$ is clearly a subgroup of $G$.  Let $L \supset K$ be a finite extension with $H=\Gal(\Kbar/L)$.  Then by Remark \ref{rmkflddef}, the field $L$ is a field of definition for $S$ if and only if the following statement holds:
\begin{center} 
For all $\sigma \in G$ and $\tau \in S$, if $\sigma|_L=\tau|_L$ then $\sigma \in S$.
\end{center}
Note $\sigma|_L=\tau|_L$ if and only if $\tau^{-1}\sigma \in H$, therefore $L$ is a field of definition if and only if
for all $\tau \in S$, we have $\tau H \subset S$, which holds if and only if $S H=S$, i.e., $H \subset H(S)$, or equivalently $L \supset \Kbar^{H(S)}=L(S)$.  Since a field of definition for $S$ exists by Lemma \ref{flddefexists}, we see that $L(S)$ is a finite extension of $K$.  Therefore $L(S)$ is the minimal field of definition for $S$.
\end{proof}

We now relate this notion to representation of primes.  Let $K$ be a number field.  Let $\Pi$ be the set of equivalence classes of sets of primes of $K$, where two sets are equivalent if they differ only by a finite set.  To every open and closed set $S \subset G$ which is closed under conjugation, we can associate a set $\calP(S)$ of primes of $K$: namely, if $L$ is a field of definition for $S$, we associate the set
\[ \calP(S)=\{\frakp\text{ a prime of $K$}:\frakp\nmid \disc(L/K),\ \Frob_\frakp \subset S|_{L}\}, \]
where $\Frob_\frakp$ is the Frobenius class at the prime $\frakp$.  If $M$ is another field of definition for $S$, then the two sets given by $L$ and $M$ differ by only a finite set, contained in the set of primes that ramify in $L$ or in $M$, and hence we have a well-defined element $\calP(S) \in \Pi$.

\begin{lem} \label{Pinj}
The above association $S \mapsto \calP(S)$ is injective.  The minimal field of definition for $S$ is Galois over $K$.
\end{lem}

\begin{proof}
Suppose that $S \neq S'$.  By Remark \ref{rmkflddef}, the compositum of a field of definition for $S$ and for $S'$ is a field of definition for both.  Therefore there exists a common field of definition $L$ for $S,S'$ which by the same remark we may take to be Galois over $K$, hence $S|_L \neq S'|_L$.  Suppose then that $\sigma \in S|_L \setminus S'|_L$; by the Chebotarev density theorem \cite[p.~169]{Lang}, there exist infinitely many primes $\frakp$ of $K$ such that $\Frob_\frakp$ is equal to the conjugacy class of $\sigma$, which is disjoint from $S'|_L$ since $S'$ is closed under conjugation.  Therefore $\calP(S) \neq \calP(S')$.

For the second statement, let $S$ be a set with minimal field of definition $L$ and let $\alpha \in G$.  Then the set $\alpha S \alpha^{-1}$ has minimal field of definition $\alpha L$: we have $\alpha\sigma\alpha^{-1}|_{\alpha L} = \alpha\tau\alpha^{-1}|_{\alpha L}$ if and only if $\sigma|_L=\tau|_L$.  Therefore if $S$ is closed under conjugation then $\alpha L=L$ and the minimal field of definition is Galois over $K$.
\end{proof}

\section{Characterizing equivalence via class groups}

In this section, we characterize the class groups which can arise from a pair of quadratic forms which represent almost the same primes.  In particular (Proposition \ref{thmclgrp}), if the forms have different fundamental discriminants, we show that they must either be of exponent dividing $2$ or of type $(2,\dots,2,4)$.  This proposition allows us to give necessary and sufficient conditions for the existence of such pairs with different fundamental discriminants (Theorem \ref{corclgrp}) and the same fundamental discriminant (Proposition \ref{iffsamedisc}).

Throughout the following sections, we will utilize the following notation.

\begin{notn} \label{qnotn}
Let $Q$ denote a (primitive, $GL_2(\bfZ)$-reduced, integral positive definite binary quadratic) form of discriminant $D=df^2$, where $d<0$ is the fundamental discriminant.  Let $K=\bfQ(\sqrt{D})=\bfQ(\sqrt{d})$, and let $R$ be the ring class field of $K$ of modulus $f$ with $h(D)=\#\Cl(D)=[R:K]$ and genus class field $P \supset K$.  By Proposition \ref{ringclassfield}, the form $Q$ corresponds to an ideal class $[\fraka]$ and to an element $\sigma \in \Gal(R/K)$.  We define the set
\[ S=\res_{R}^{-1}(\{\sigma,\sigma^{-1}\}) \subset \Gal(\bfQbar/\bfQ). \]
Note that $\calP(S)$ (as in Lemma \ref{Pinj}) is the set of primes represented by $Q$, up to a finite set (contained in the set of primes dividing $f$).  
\end{notn}

The set $S$ is open and closed in $\Gal(\bfQbar/\bfQ)$ and closed under conjugation.  Let $L=L(S)$ be the minimal field of definition for $S$, which exists by Corollary \ref{minflddefexists}; since $R$ is a field of definition for $S$, we have $L \subset R$.  (Note here we take the base field in \S 2 to be $\bfQ$.)

\begin{lem} \label{lem24}
We have $[R:L] \leq 2$, and $[R:L]=2$ if and only if $\sigma|_L$ has order $2$ and $\sigma$ has order $4$.  Moreover, we have $P \subset L$.
\end{lem}

\begin{proof}
Since $S|_R=\{\sigma,\sigma^{-1}\}$, we have
\[ 2 \geq \# S|_R=[R:L](\# S|_L), \]
so $[R:L] \leq 2$.  Moreover, $[R:L]=2$ if and only if $\#S|_R=2$ and $\#S|_L=1$, which holds if and only if $\sigma|_L=\sigma^{-1}|_L$ and $\sigma \neq \sigma^{-1}$, i.e., $\sigma|_L$ has order $2$ and $\sigma$ has order $4$.

To prove that $P \subset L$, note that in either case $\Gal(R/L)$ is generated by $\sigma^2 \in \Cl_f(d)^2=\Gal(R/P)$.
\end{proof}

Now suppose that $Q_1$ and $Q_2$ are a pair of forms, following Notation \ref{qnotn} with appropriate subscripts.  It is immediate from Lemma \ref{repsame} that $Q_1$ and $Q_2$ have the same set $\calP(S)$ (up to a finite set) and by the injectivity of Lemma \ref{Pinj} the same set $S$, hence the same minimal field of definition $L$.

\begin{lem} \label{samegenus}
If $Q_1 \sim Q_2$, then we have $K_1K_2 \subset L$, and $K_1K_2$ is fixed by all elements of $S$.  Moreover, we have equality of genus class fields $P_1=P_2$.
\end{lem}

\begin{proof}
This follows immediately from the fact that $K_i \subset P_i \subset L$ and that $P_i$ is the maximal subextension of $L/\bfQ$ of exponent dividing $2$.  
\end{proof}

We denote this common genus class field by $P=P_1=P_2$.

\begin{cor} \label{samegenuscor}
If $Q_1 \sim Q_2$, then $\sigma_1|_P=\sigma_2|_P$.
\end{cor}

\begin{proof}
Note $\sigma_2|_P=\sigma_2^{-1}|_P$.  Since $P \subset L$, by Lemma \ref{repsame} we conclude $\sigma_1|_P = \sigma_2|_P$.
\end{proof}

We now distinguish two cases, depending on whether $Q_1,Q_2$ have the same fundamental discriminant.

\begin{prop} \label{thmclgrp}
Suppose $Q_1 \sim Q_2$ and $K_1 \neq K_2$.  Then for $i=1,2$, the group\/ $\Gal(R_i/K_i)$ is of type dividing $(2,\dots,2,4)$, and the minimal field of definition is equal to the common genus class field, i.e.~$L=P$.
\end{prop}

\begin{proof}
Let $\alpha \in \Gal(L/\bfQ)$ be any element of order not dividing $2$.  From Proposition \ref{ringclassfield} we have
\[ \Gal(L/\bfQ)=\Gal(L/K_i) \rtimes \Gal(K_i/\bfQ) \]
where the nontrivial element of $\Gal(K_i/\bfQ) \cong \bfZ/2\bfZ$ acts on $\Gal(L/K_i)$ by inversion.  Suppose that $\alpha \in \Gal(L/\bfQ)$ is an element of order $>2$.  Then in fact $\alpha \in \Gal(L/K_i)$, since every element of $\Gal(L/\bfQ) \setminus \Gal(L/K_i)$ has order $2$.  Therefore the centralizer of $\alpha$ in $\Gal(L/\bfQ)$ is the group $\Gal(L/K_i)$.  Hence if such an $\alpha$ exists, then $K_i$ is determined by $L$, so $K_1=K_2$.  So $K_1 \neq K_2$ implies that $\Gal(L/\bfQ)$ is of exponent $2$, and then from the exact sequence
\[ 0 \to \Gal(R_i/L) \to \Gal(R_i/K_i) \to \Gal(L/K_i) \to 0 \]
and the fact that $[R_i:L] \leq 2$ we see that $\Gal(R_i/K_i)$ is of type dividing $(2,\dots,2,4)$.

The second statement then follows, since then $L \subset P$.  
\end{proof}

\begin{rmk}
This proposition answers a question of Jagy and Kaplansky \cite{JagyKaplansky}.  Two ideal classes are said to be in the same \emph{genus} if their ratio is a square of an ideal class.  Jagy and Kaplansky call a form $Q$ \emph{bi-idoneal} if its genus consists of only $Q$ and its inverse; in their terminology, every ``non-trivial'' pair of forms (i.e., $d_1 \neq d_2$) representing the same primes they found was bi-idoneal.  

Proposition \ref{thmclgrp} shows that this always holds: if $Q_1,Q_2$ represent the same primes outside a finite set and $d_1 \neq d_2$, then $Q_1$ and $Q_2$ are bi-idoneal.  This follows from the fact that a finite abelian group $G$ has $\#(G^2) \leq 2$ if and only if $G$ is of type dividing $(2,\dots,2,4)$.
\end{rmk}

We can now formulate necessary and sufficient conditions for the existence of pairs which represent almost the same primes with different fundamental discriminants.

\begin{thm} \label{corclgrp}
Let $Q_1,Q_2$ be forms, and suppose that $K_1 \neq K_2$.  Then $Q_1 \sim Q_2$ if and only if both of the following hold:
\begin{enumroman}
\item $R_1$ and $R_2$ have the same genus class field $P$, and
\[ \sigma_1|_P=\sigma_2|_P \in \Gal(P/K_1K_2); \]
\item For $i=1,2$, the group $\Gal(R_i/K_i)$ is either of exponent dividing $2$, or is of type $(2,\dots,2,4)$ and $\sigma_i$ has order $4$.
\end{enumroman}
\end{thm}

\begin{proof}
We have shown these conditions are necessary: condition (i) follows from Lemma \ref{samegenus} and Corollary \ref{samegenuscor} and (ii) follows from Proposition \ref{thmclgrp}.

Now we show that these conditions are also sufficient.  For $i=1,2$, let $L_i$ be the minimal field of definition of $S_i$ (as in Notation \ref{qnotn}, with subscripts).  From Lemma \ref{lem24} and (i), we have $P \subset L_i$, and since $R_i$ is a field of definition for $S_i$ we have $L_i \subset R_i$.  We now will show that in fact $L_i=P$.  From (ii), either $\Gal(R_i/K_i)$ is of exponent dividing $2$ and $R_i=L_i=P$ already, or $\Gal(R_i/K_i)$ is of type $(2,\dots,2,4)$ and $\sigma_i$ has order $4$.  But then $P$ is a field of definition for $S_i$, since $\res_R^{-1}(\sigma_i|_P)=\{\sigma_i,\sigma_i^{-1}\}$, hence $L_i \subset P$, so $L_i=P$ in this case as well.  Therefore $L_1=L_2=L$.  

Now let $p$ be a prime which is unramified in $R_1R_2$.  Then $\sigma_1 \in \Frob_p|_L$ if and only if $\sigma_2 \in \Frob_p|_L$, so then $Q_1 \sim Q_2$ by Lemma \ref{repsame}.
\end{proof}

To conclude this section, we consider the case when two forms have the same fundamental discriminant.

\begin{prop} \label{iffsamedisc}
Let $Q_1,Q_2$ be forms with $d_1=d_2=d$.  

Suppose that $f_1 \mid f_2$, and let 
\[ \phi:\Cl(D_2) \to \Cl(D_1) \]
be the natural (restriction) map.  Then $Q_1 \sim Q_2$ if and only if $\phi(\sigma_2)=\sigma_1$ and one of the following holds: either $\phi$ is an isomorphism, or
\[
\text{The kernel of $\phi$ has order $2$, generated by $\sigma_2^2$, and $\sigma_1$ has order $2$.}
\tag{$\dagger$} \]

More generally, we have $Q_1 \sim Q_2$ if and only if there exists a form $Q$ of discriminant $D=df^2$ with $Q_1 \sim Q \sim Q_2$, where $f=\gcd(f_1,f_2)$.  
\end{prop}

\begin{proof}
From Proposition \ref{ringclassfield}, we conclude that $R_1 \subset R_2$.  If $R_1=R_2$ then $\phi$ is an isomorphism.  Otherwise, by Lemma \ref{lem24}, we have $[R_2:R_1]=2$ and $Q_1 \sim Q_2$ if and only if $\res_{R_2}^{-1}(\sigma_1)=\{\sigma_2,\sigma_2^{-1}\}$, where $\sigma_2$ has order $4$ and $\sigma_1$ has order $2$.  Now $\sigma_1$ has order $2$ if and only if $\sigma_2^2 \in \ker \phi$, and $\ker \phi$ is generated by $\sigma_2^2$ if and only if $\sigma_2$ has order $4$, which is condition $(\dagger)$.  This proves the first statement.

To prove the second statement, let $R=R_f$.  Then by Corollary \ref{rgcd}, $R_1 \cap R_2 = R$.  Since $L \subset R_1,R_2$ we see that $L \subset R$, therefore by Remark \ref{rmkflddef} the field $R$ is a field of definition for $S$.  Let $Q$ be the form of discriminant $df^2$ associated to $\sigma_1|_R$.  Again by Lemma \ref{lem24}, we see that either $R_1=R$, in which case $Q_1 \sim Q$, or $[R_1:R]=2$, in which case $L=R$ and as above we have $Q_1 \sim Q$.  Similarly, let $Q'$ be the form of discriminant $df^2$ associated to $\sigma_2|_R$.  Then $Q_2 \sim Q'$.  Since $Q_1 \sim Q_2$, we have $Q \sim Q'$.  But $Q$ and $Q'$ have the same discriminant, which implies that $Q=Q'$, by Remark \ref{D1D2}.
\end{proof}

\section{Bounding class groups}

Recall as in the introduction, to every equivalence class $C$ of forms, we associate the set $\delta(C)$ of fundamental discriminants of the forms in $C$ as well as the set $\Delta(C)$ of discriminants of forms in $C$.  In this section, we will prove that there are only finitely many equivalence classes $C$ with $\#\delta(C) \geq 2$.  More precisely, we will prove the following statement.  

\begin{prop} \label{d1neqd2discs}
The sets
\[ \calD_\delta=\bigcup_{\#\delta(C) \geq 2} \delta(C) \quad \text{ and } \quad \calD_\Delta=\bigcup_{\#\delta(C) \geq 2} \Delta(C), \]
are finite and effectively computable.  Moreover, $\# \calD_\delta \leq 226$ and $\# \calD_\Delta \leq 425$.
\end{prop}

First note the following lemma.

\begin{lem}[{\cite[Lemma 5]{Weinberger}}] \label{classbound}
Let $K=\bfQ(\sqrt{d})$ have discriminant $d<0$, let $\fraka$ be an integral ideal of $K=\bfQ(\sqrt{d})$ and let $c$ be a positive integer such that $\fraka^c$ is principal.  If $\fraka$ is not a principal ideal generated by a rational integer and $\fraka$ is prime to $d$, then $(N\fraka)^c > |d|/4$.
\end{lem}

To prove this lemma, one shows that if $(\alpha)=\fraka^c$, then $\alpha$ is not a rational integer by considering the factorization of $\fraka$ in $K$, and therefore $N(\fraka^c)=N(\alpha)^c > |d|/4$.  

\begin{cor} \label{expc}
If $\Cl(d)$ has exponent $c$, then for all primes $p$ such that $p^c \leq d/4$ we have $(d/p) \neq 1$.
\end{cor}

\begin{proof}
Suppose that $(d/p)=1$; then $(p)=\frakp \overline{\frakp}$ in the ring of integers $A$ of $K=\bfQ(\sqrt{d})$.  Since $N\frakp=p$ is not a square, we know that $\frakp$ is not generated by a rational integer.  The lemma implies then that $(N\frakp)^c=p^c > d/4$.
\end{proof}

\begin{lem} \label{2224restrictsf}
If $\Cl_f(d)$ is of type dividing $(2,\dots,2,4)$ and $|d|>2500$, then $f \in \{1,2,3,4,6,8,12\}$.
\end{lem}

\begin{proof}
Recall the exact sequence of Lemma \ref{exactseq}
\[ 1 \to (A/fA)^*/(\bfZ/f\bfZ)^* \to \Cl_f(d) \to \Cl(d) \to 1, \]
where note that $|d|>4$ implies $A_f^*=A^*$.  

Since the map $\Cl_f(d) \to \Cl(d)$ is surjective, we see that $\Cl(d)$ is itself of type dividing $(2,\dots,2,4)$.  Let $p$ be an odd prime such that $p \mid f$.  From Proposition \ref{CRT}, we conclude that $p^2 \nmid f$ and $p=3$ or $p=5$.  When $|d|>2500$, or equivalently when $|d/4|>5^4$, we cannot have $5 \mid f$, for this can happen only if $(d_i/5)=1$, which contradicts Corollary \ref{expc}.  If $2 \mid f$, then since $(d/2)=1$ cannot occur, and $(d/2)=-1$ implies $3 \mid \Cl_f(d)$, we must have $(d/2)=0$.  But then again from the proposition we see that $16 \nmid f$ and $24 \nmid f$.  
\end{proof}

Let $Q_1,Q_2$ be forms with $d_1 \neq d_2$.  Let $K_0$ be the real quadratic field contained in $K_1K_2$.  


\begin{lem} \label{inertyeah}
Let $Q_1 \sim Q_2$ and suppose $|d\sbmin|=\min\{|d_1|,|d_2|\}>2500$.  Then 
\[ K_0 \in \{\bfQ(\sqrt{2}),\bfQ(\sqrt{3}),\bfQ(\sqrt{6})\}.  \]
Moreover, if $p^4 \leq |d\sbmin|/4$ and $p$ is inert in $K_0$, then $p$ ramifies in $K_1$ and $K_2$.
\end{lem}

\begin{proof}
By Lemma \ref{samegenus}, the ring class fields $R_1$ and $R_2$ have the same genus class field, and by Lemma \ref{thmclgrp}, the group $\Cl_{f_i}(d_i)$ is of type dividing $(2,\dots,2,4)$ for $i=1,2$.  By Corollary \ref{samecor}, the same set of odd primes divide the discriminants $D_1,D_2$.  Then by Lemma \ref{2224restrictsf}, we see that $d_1/d_2 \in 2^{\bfZ} 3^{\bfZ}$.  Therefore the discriminant of $K_0$ is supported only at the primes $2$ and $3$, and $K_0$ is one of the fields listed.

Let $p$ be a prime with $p^4 \leq d\sbmin/4$ which is inert in $K_0$.  We know that $(d_1/p),(d_2/p) \neq 1$, by Corollary \ref{expc}.  We cannot have $(d_1/p)=(d_2/p)=-1$, as then $(d_1d_2/p)=1$ so $p$ splits in $K_0$.  Therefore say $(d_1/p)=0$; then $p$ is ramified in $K_1$ so $p$ is ramified in $K_1K_2=K_0K_2$, so $p$ is ramified in $K_2$ as well.
\end{proof}

\begin{rmk} \label{inertyeahrmk}
This lemma proves that given a fundamental discriminant $d$ with $|d|>2500$, one can explicitly determine all possibilities for fundamental discriminants $d'$ of forms $Q'$ with $Q' \sim Q$.
\end{rmk}

\begin{lem} \label{cntm}
Let $p_1=3$, $p_2=5$, $\dots$ be the sequence of odd primes in increasing order, and for each $t \in \bfZ_{\geq 1}$ let
\[ \dhat_t=4p_1 \dots p_{t-1}. \]
Let $d < -3$ be a fundamental discriminant with $g$ distinct prime factors, and let $t \in \bfZ_{\geq 1}$.  Then
\[ |d| \geq \dhat_t p_t^{g-t}. \] 
\end{lem}

\begin{proof}
First, we prove that $|d| \geq \dhat_g$.  If $d \equiv 0 \pmod{4}$, then this is clear.  If $d \equiv 1 \pmod{4}$ and $g=1$, then by assumption $|d| \geq 7>4$.  If $g \geq 2$, then $p_g \geq 5$, and therefore
\[ |d| \geq p_1 \dots p_g \geq 4p_1 \dots p_{g-1}. \]
It then follows that $|d| \geq \dhat_g \geq \dhat_t p_t^{g-t}$ for $g \geq t$.  But for $g<t$, we also have
\[ |d| \geq \dhat_g=\frac{\dhat_t}{p_{g+1} \cdots p_t} \geq \frac{\dhat_t}{p_t^{t-g}} \]
as claimed.
\end{proof}

By the preceding two lemmas, we can apply the result of Tatuzawa (Lemma \ref{tatuzawa}) to obtain the following.

\begin{prop} \label{bndtatusol2}
Let $Q_1,Q_2$ be forms representing almost the same primes such that $d_1 \neq d_2$.  Then we have $\min\{|d_1|,|d_2|\} \leq B=80604484=4\cdot 67^4$.
\end{prop}

\begin{proof}
Apply Lemma \ref{tatuzawa} with $\epsilon=1/\log B$.  Note that $\log B >11.2$.  Since there is at most one possible exceptional discriminant, we may assume without loss of generality that $d=d_1$ is not exceptional, hence
\[ h(d) > \left(\frac{0.655}{\pi}\right)\frac{|d|^{1/2-1/\log B}}{\log B}. \]

We suppose that $|d| > B$ and derive a contradiction.  By Lemma \ref{inertyeah}, every prime $p \leq 67$ which is inert in $K_0$ must divide $d$.  Let $g$ be the number of distinct prime factors of $d$; since $\#\Cl(d)[2]=2^{g-1}$ (Corollary \ref{2gm1}) and $\Cl(d)$ is of type dividing $(2,\dots,2,4)$, we see that $h(d) \leq 2^g$.  

For $b \in \bfZ_{>0}$, let
\[ d_0(b,q) = \prod_{\substack{2 < p \leq b \\ (p/q)=-1}} \negthickspace\negthickspace p. \]
From Lemma \ref{inertyeah}, we have three cases to consider.  If $K_0=\bfQ(\sqrt{2})$, then $p$ is inert in $K_0$ if and only if $p \equiv 3,5 \pmod{8}$.  Therefore
\[ d_0(67,8)=\prod_{\substack{p \leq 67 \\ p \equiv 3,\,5 \psod{8}}} \negthickspace\negthickspace p = 3 \cdot 5 \cdot \ldots \cdot 61 \cdot 67 > 2.4 \cdot 10^{16}, \] 
and by Lemma \ref{inertyeah}, we have $d_0(67,8) \mid d$, so $|d| \geq d_0(67,8)$.  For $K_0=\bfQ(\sqrt{3})$, the prime $p$ is inert in $K_0$ if and only if $p \equiv 5,7 \pmod{12}$, so $d_0(67,12)=5 \cdot \ldots \cdot 53 \cdot 67$, and $d_0(67,12) \mid d$ so $|d| > 6.3 \cdot 10^{13}$.  In a similar way, for $K_0=\bfQ(\sqrt{6})$, we obtain $d_0(67,24)=7 \cdot 11 \cdot \ldots \cdot 61 > 2.8 \cdot 10^{13}$.  

In any case, we see that $|d| > 2.8 \cdot 10^{13}$, and hence
\[ 2^g \geq h(d) > \left(\frac{0.655}{\pi \log B}\right)(2.8 \cdot 10^{13})^{1/2-1/\log B} > 10897 \]
so $g \geq 14$.

By Lemma \ref{cntm}, we have $|d| \geq \dhat_{14} \cdot 47^{g-14}$, where $\dhat_{14} > 2.6 \cdot 10^{17}$.  But this implies that
\begin{align*}
2^g \geq h(d) &> \left(\frac{0.655}{\pi \log B}\right)\dhat_{14}^{1/2-1/\log B} \left(47^{1/2-1/\log B}\right)^{g-14} \\
&> 226989 \cdot 2^{g-14},
\end{align*}
which is a contradiction.
\end{proof}

\begin{proof}[Proof of Proposition ${\ref{d1neqd2discs}}$]
First, by an exhaustive list, we find that there are exactly $226$ fundamental discriminants $d$ with $|d| \leq B$ such that $\Cl(d)$ is of type dividing $(2,\dots,2,4)$.  To speed up this computation, we use Corollary \ref{expc} to rule out many of these discriminants.  This was accomplished in \textsf{MAGMA}.  (The code is available from the author by request.)  By Proposition \ref{bndtatusol2}, we have missed at most one possible fundamental discriminant from the set $\calD_d$.

Next, we show that there are exactly $199$ discriminants $D=df^2$ of nonmaximal orders with $|d| \leq B$ such that $\Cl_f(d)$ is of type dividing $(2,\dots,2,4)$.  By Lemma \ref{2224restrictsf}, we know that $f \in \{2,3,4,6,12\}$.  We can use any algorithm which computes class groups (e.g. enumeration) to check these finitely many nonmaximal orders.

Now suppose that $Q_1,Q_2$ are forms that represent the same primes with $|d_1| < |d_2|$.  Then $|d_1| \leq B$, and we must show that $|d_2| \leq B$ as well to have computed $\calD_d$ and therefore $\calD_D$ as well.  If $|d_1| \leq 2500$, then from the list of discriminants we see that $|D_1| \leq 29568$; since $\bfQ(\sqrt{d_2}) \subset P_1$, we see from Lemma \ref{samegenus} that $|d_2| \leq 4\cdot 29568 < B$.  Otherwise, by Remark \ref{inertyeahrmk}, there are only $3$ possibilities for $d_2$, and since $|d_1| \leq 10920$, it follows that $|d_2| \leq 12\cdot 10920 \leq B$ as well, completing the proof.
\end{proof}

\section{Finding the pairs of quadratic forms}

To conclude, we list all forms with $K_1 \neq K_2$.  Using Corollary \ref{computegenus}, we first compute the genus class field for each of the $425$ discriminants found in section 5.  We find $86$ pairs of discriminants for which the genus class fields are equal.

We now apply Theorem \ref{corclgrp}.  If the class group of both discriminants are both of exponent 2, then for every $\sigma \in \Gal(R/K_1K_2)$, we obtain a pair corresponding to $\sigma_i=\sigma \in \Gal(R/K_i)$.  For each $i$ such that $\Gal(R_i/K_i)$ has a factor $\bfZ/4\bfZ$, we proceed as follows: for each $\sigma \in \Gal(R_i/K_1K_2) \subset \Gal(R_i/K_i)$ of order $4$, we compute the fixed field of $\sigma|_P$ by finding a prime $p \nmid D_i$ represented by the form $Q \leftrightarrow \sigma$, and compute (using Legendre symbols) the largest subfield of $P$ in which $p$ splits completely.  Then every pair $\sigma_1,\sigma_2$ with the same fixed subfield (so that $\sigma_1|_P=\sigma_2|_P$) gives rise to a pair of forms.

\begin{exm}
The discriminants $D_1=-1056=-264\cdot 2^2$ and $D_2=-2112=-132\cdot 4^2$ give rise to the common genus class field $P=\bfQ(i,\sqrt{2},\sqrt{-3},\sqrt{-11})$ each with class group of type $(2,2,4)$.  The forms of order $4$ of discriminant $-1056$ are
\[ \la 5, 2, 53 \ra, \la 15, 12, 20 \ra, \la 7, 6, 39 \ra, \la 13, 6, 21 \ra, \]
and those of discriminant $-2112$ are
\[ \la 17, 8, 32 \ra, \la 21, 18, 29 \ra, \la 7, 4, 76 \ra, \la 19, 4, 28 \ra. \]
The first form $\la 5, 2, 53 \ra$ represents the prime $5$, so we compute the Legendre symbols 
\[ (-1/5),(2/5),(-3/5),(-11/5), \]
and find the fixed field $\bfQ(i,\sqrt{6},\sqrt{-11}) \subset P$.  Continuing in this way, we find that only the pair $\la 7,6,39 \ra$ and $\la 7,4,76 \ra$ have a common fixed field, namely the field $\bfQ(\sqrt{2},\sqrt{-3},\sqrt{11})$, and this proves that they represent the same primes (those which are congruent to $7,79,127,151,175 \pmod{264}$).
\end{exm}

Carrying out this calculation for each of the $86$ pairs, and supplementing this list with any pairs arising from the same fundamental discriminant, we obtain the forms listed in Tables 1--5.

\begin{thm} \label{thm1}
There are exactly\/ $67$ equivalence classes of forms $C$ such that $\# \delta(C) \geq 2$.  There are exactly\/ $6$ classes with $\#\delta(C)=3$ and there is no class with $\# \delta(C) \geq 4$.  
\end{thm}

\begin{defn} \label{except}
The \emph{exceptional set} $E$ of a form $Q$ is the set of primes $p$ such that $Q$ represents $p$ and there exists a form $Q' \sim Q$ such that $Q'$ does not represent $p$.
\end{defn}

\begin{rmk} \label{jkmiss}
Jagy and Kaplansky \cite{JagyKaplansky} miss the two pairs 
\[ \la 5, 0, 6 \ra, \la 11, 4, 14 \ra \quad \text{ and } \quad \la 3, 0, 40 \ra, \la 27, 12, 28 \ra \] 
in their ``near misses'' (those pairs with exceptional set not contained in $\{2\}$).  Moreover, the form $\la 4,4,9 \ra$ in their paper should be $\la 4,4,19 \ra$.
\end{rmk}

\section{Forms with the same fundamental discriminant}

In this section, we treat the case when the forms have the same fundamental discriminant.  We will again use Notation \ref{qnotn}.  Throughout, let $Q_1,Q_2$ be forms with $d_1=d_2=d<0$.  

If $f_1=f_2$, so that $D_1=D_2$, then by Remark \ref{D1D2} either $Q_1=Q_2$ or $Q_1 \not\sim Q_2$.  So without loss of generality we may assume that $f_1 < f_2$.  

We begin with a general lemma about quadratic forms.

\begin{defn} \label{rlift}
Let $Q$ be a form of discriminant $D<0$ and let $r \in \bfZ_{\geq 1}$.  The form $Q'$ is an \emph{$r$-lift} of $Q$ if the following conditions hold:
\begin{enumalph}
\item $Q$ and $Q'$ have the same fundamental discriminant $d=d'$;
\item The discriminant of $Q'$ satisfies $D'=r^2 D$; 
\item In the natural (restriction) map
\[ \phi:\Cl(D') \to \Cl(D) \]
we have $\phi(\sigma')=\sigma$, where $\sigma \leftrightarrow Q$ and $\sigma' \leftrightarrow Q'$.  
\end{enumalph}
\end{defn}

\begin{lem} \label{char2powerforms}
Let $Q=\la a,b,c \ra$ be an $SL_2(\bfZ)$-reduced form associated to $\sigma$.  Then $\sigma$ has order dividing $2$ if and only if $0=b$ or $b=a$ or $a=c$.

Suppose that $\sigma$ has order dividing $2$ and $2 \mid D$.  Then $Q$ has a $2$-lift $Q'$ with $Q' \leftrightarrow \sigma'$ of order $2$ if and only if $0=b$.
\end{lem}

\begin{proof}
Throughout this proof, we require only that forms be $SL_2(\bfZ)$-reduced rather than $GL_2(\bfZ)$-reduced, but we maintain all other assumptions on our forms, as in the introduction.  Recall that $Q$ is $SL_2(\bfZ)$-reduced if and only if $|b| \leq a \leq c$ and $b=0$ if either $|b|=a$ or $a=c$.

The first statement of the lemma is classical: The opposite of the form $Q$ is the form $SL_2(\bfZ)$-equivalent to $Q'=\la a,-b,c \ra$.  But this form is already $SL_2(\bfZ)$-reduced, unless $|b|=a$ or $a=c$, and in either of these cases in fact $Q'$ is $SL_2(\bfZ)$-equivalent to $Q$, so that $\sigma$ has order dividing $2$.

For the second statement, first suppose $0=b$ and that $a$ is odd.  Note that the form $Q'=\la a,0,4c \ra$ is a $2$-lift of $Q$, since the set of primes which it represents is a subset of those represented by $Q$.  If $c$ is odd, then a $2$-lift is $\la 4a,0,c \ra$ if $4a \leq c$ and $\la c,0,4a \ra$ if $4a>c$.  This concludes this case, because if $a$ and $c$ are both even then $Q$ is not primitive.  

Next, suppose that $b=a$.  Then since $D$ is even, $a$ is even, so $c$ is odd.  Therefore a $2$-lift of $Q$ is the form $SL_2(\bfZ)$-equivalent to $Q'=\la 4a,2a,c \ra$, which is either $Q'$ if $4a<c$, or $\la c,-2a,4a\ra$ if $2a<c<4a$, or $\la c,2(c-a),4a+c \ra$ if $c<2a$; we cannot have $4a=c$ or $2a=c$ as then $c$ is even and $Q$ is not primitive.  In any case, the $2$-lift visibly has order $>2$, therefore all $2$-lifts have order $>2$ since they differ by an element of the kernel which is of order dividing $2$, by Proposition \ref{CRT}.

Finally, suppose $a=c$.  Here, we know that $b$ is even so $a$ is odd, and a $2$-lift of $Q$ is the form $SL_2(\bfZ)$-equivalent to $Q'=\la a,2b,4a \ra$, which is $Q'$ if $2b<a$ and $\la a,2(b-a),5a-2b \ra$ if $2b>a$; we cannot have $2b=a$, since $a$ is odd.  This form has order dividing $2$ if and only if $b=a$ which is impossible ($a$ must be even from the previous paragraph), and otherwise this lift has order $>2$.
\end{proof}

\begin{prop} \label{unique2or4}
Let $Q_1,Q_2$ be forms with $d_1=d_2=d$ and $f_1<f_2$.  Then $Q_1 \sim Q_2$ if and only if $Q_2$ is the unique $2$- or the unique $4$-lift of $Q_1$.  
\end{prop}

\begin{proof}
First, suppose that $f_1 \mid f_2$ and that $A_{f_1}$ and $A_{f_2}$ have the same number of roots of unity.  Note that the set of primes represented by $Q_2$ is contained in the set of primes represented by $Q_1$ up to a finite set if and only if $Q_2$ is an $r$-lift of $Q_1$ for some $r \in \bfZ_{>1}$.  Moreover, if there exist two such $r$-lifts $Q_2,Q_2'$, then these two forms will represent disjoint, infinite nonempty sets of primes.  Putting these together, we see that $Q_1 \sim Q_2$ if and only if $Q_2$ is the unique $r$-lift of $Q_1$ for some $r \in \bfZ_{>1}$.

From Lemma \ref{lem24} we have $[R_2:R_1] \in \{1,2\}$.  On the other hand, by Lemma \ref{samegenus}, $R_1$ and $R_2$ have the same genus class field, so from Proposition \ref{computegenus}, if $p \nmid d$ is an odd prime then $p \mid f_1$ if and only if $p \mid f_2$.  From Lemma \ref{exactseq} we have
\begin{equation}
[R_2:R_1]=\frac{h(D_2)}{h(D_1)}=u\frac{f_2}{f_1} \in \{1,2\} \tag*{$(*)$}
\end{equation}
where
\[ u=
\begin{cases}
\displaystyle{\left(1-\legen{d}{2}\frac{1}{2}\right)}, & \text{ if } 2 \nmid f_1\text{ and }2 \mid f_2; \\
1, & \text{ otherwise}.
\end{cases} \]

From Proposition \ref{iffsamedisc}, there exists a form $Q$ of discriminant $df^2$ with $f=\gcd(f_1,f_2)$ such that $Q_1 \sim Q \sim Q_2$.  But since $u \in \frac{1}{2}\bfZ$ we see from $(*)$ that $f_i/f \in 2^{\bfZ}$ for $i=1,2$, so $f_2/f_1 \in 2^{\bfZ}$ as well and hence since $f_1<f_2$ we have $f=f_1 \mid f_2$ and $Q=Q_1$.  Moreover, we have $u=1/2$ or $u=1$ and hence either $f_2=2f_1$ or $f_2=4f_1$, so $Q_2$ is the unique $2$- or $4$-lift of $Q_1$.

To conclude, suppose that the two orders have different numbers of roots of unity.  Then $d=-3,-4$ and $A_{f_1}$ is the maximal order and $A_{f_2}$ is not.  Repeating the above analysis, we easily verify that either $f_2=2f_1$ or $f_2=4f_1$; the finitely many cases that can occur are listed in Table 6.
\end{proof}

To conclude, from this proposition it suffices to give necessary and sufficient conditions for the form $Q_1$ to have a unique $2$- or $4$-lift.  Note that if $Q_2$ is the unique $4$-lift of $Q_1$, and $Q$ is the unique $2$-lift of $Q_1$, then in fact $Q_2$ is the unique $2$-lift of $Q$, and $Q_1 \sim Q \sim Q_2$.  Therefore it suffices to give criteria equivalent to those occurring in Proposition \ref{iffsamedisc}.

\begin{thm} \label{k1eqk2}
Let $Q_1=\la a_1,b_1,c_1 \ra$ be a form.  Then there exists a form $Q_2 \sim Q_1$ such that $|D_2|>|D_1|$ and $d_2=d_1=d$ if and only if one of the following holds:
\begin{enumroman}
\item $d \equiv 1 \pmod{8}$ and $2 \nmid D_1$;
\item $2 \mid D_1$ and either $b_1=a_1$ or $a_1=c_1$;
\item $d=-3$ and $Q_1 \in \{\la 1,1,1 \ra,\la 1,0,3 \ra \}$;
\item $d=-4$ and $Q_1=\la 1,0,1 \ra$.
\end{enumroman}
\end{thm}

\begin{proof}
If $d=-3$ or $d=-4$, we refer to Proposition \ref{unique2or4} (and Table 6) and find cases (iii) and (iv).  

More generally, we apply Proposition \ref{iffsamedisc}.  The map $\phi$ is an isomorphism if and only if $h(D_2)=h(D_1)$.  By Proposition \ref{CRT}, this occurs if and only if $(d/2)=1$ (and $f_2=2f_1$), which is case (i).  

For condition $(\dagger)$ from Proposition \ref{iffsamedisc}, first for any positive integer $f$, let 
\[ C(f) = \frac{(A/fA)^*}{(\bfZ/f\bfZ)^*}. \]
From the functoriality of the exact sequence of Lemma \ref{exactseq}, we obtain a commutative diagram
\[
\xymatrix{
1 \ar[r] & C(f_2) \ar[r] \ar[d]^{\psi} & \Cl_{f_2}(d) \ar[d]^{\phi} \\
1 \ar[r] & C(f_1) \ar[r] & \Cl_{f_1}(d)
} \]
Now if $(\dagger)$ holds then $C(f_2) \to C(f_1)$ is a nonsplit $\bfZ/2\bfZ$-extension, so we see from Proposition \ref{CRT} that $2 \mid D_1$.  Therefore $(\dagger)$ holds if and only if $2 \mid D_1$, $\sigma_1$ has order $2$ and $\sigma_2$ has order $4$.  The result now follows from Lemma \ref{char2powerforms}.
\end{proof}

\section{Computing class groups}

To give an alternative proof of Proposition \ref{d1neqd2discs}, we may also characterize with at most one possible exception all imaginary quadratic extensions having class group of type dividing $(2,\dots,2,4)$.  This result is not needed in the sequel, but it also yields an independent result (Theorem \ref{alldiscrim}).

It was a classical problem to characterize field discriminants whose class group has exponent dividing $2$, comprised of quadratic forms which are said to be ``alone in their genus''.  It has long been known that the Brauer-Siegel theorem implies that there are only finitely many such discriminants \cite{Chowla}.

\begin{prop}[Weinberger \cite{Weinberger}, Louboutin \cite{Louboutin}]
The number of discriminants $D=df^2<0$ such that $\Cl_f(d)$ has exponent dividing $2$ is finite.  There are
at least $65$ and at most $66$ such fundamental discriminants, and at least $36$ and at most $37$ such discriminants of nonmaximal orders.  

Under the assumption of a suitable generalized Riemann hypothesis, there are exactly $65$ and $36$ of these discriminants, respectively.
\end{prop}

The list of these discriminants can be found in \cite[Table 5]{Borevich}.  Here we have a small variant of this problem, to which we may apply the same techniques.  

\begin{thm} \label{alldiscrim}
There are at least $226$ and at most $227$ fundamental discriminants $D=d$ such that $\Cl(d)$ is of type dividing $(2,\dots,2,4)$, and at least $199$ and at most $205$ such discriminants $D$ of nonmaximal orders.  
\end{thm}

These extensions are listed in Tables 7--16.  Our proof of the proposition will again rely on the result of Tatuzawa (Lemma \ref{tatuzawa}).

\begin{lem} \label{effconstant}
There are effectively computable constants $C_9$, $C_{10}$, and $C_{11}$ satisfying the following condition: 

With at most one exception, for all fundamental discriminants $d<0$ with $g$ distinct prime factors such that $|d| \geq C_9$ and $\Cl(d)$ is of type dividing $(2,\dots,2,4)$, we have $g \in \{10,11\}$ and $|d| <C_g$.
\end{lem}

\begin{proof}
Let $d<0$ be a fundamental discriminant with $g$ distinct prime factors and class group of type dividing $(2,\dots,2,4)$.  Recall as in the proof of Proposition \ref{bndtatusol2} that $h(d) \leq 2^g$.

Let $C_9$ be the smallest positive integer such that
\[ 2^9=512 \leq \frac{0.655}{\pi e}\frac{\sqrt{C_9}}{\log C_9} \]
(allowable, since $\sqrt{x}/\log x$ is increasing for $x \geq e^2$).  A calculation shows that $\log C_9 > 23$. Now apply Lemma \ref{tatuzawa} with $\epsilon=1/\log C_9$.  

Suppose that $d$ is not the exceptional discriminant.  Then if $|d| \geq C_9$, we have
\[ 2^g \geq h(d) > \left(\frac{0.655}{\pi}\right)\frac{|d|^{1/2-1/\log C_9}}{\log C_9}. \]
In particular, this implies that
\[ 2^g > \frac{0.655}{\pi e}\frac{\sqrt{C_9}}{\log C_9} \geq 2^9 \]
and therefore $g>9$.  

By Lemma \ref{cntm}, we have $|d| \geq \dhat_9\cdot 29^{g-9}$ and hence
\[ 2^g \geq h(d) > \left(\frac{0.655}{\pi}\right)\frac{\dhat_9^{1/2-1/\log C_9}}{\log C_9} \left(29^{1/2-1/\log C_9}\right)^{g-9}. \]
This inequality implies that $g < 12$.  

For $t \in \{10,11\}$, let $C_t$ the smallest positive integer such that
\[ 2^t \leq \left(\frac{0.655}{\pi}\right)\frac{C_t^{1/2-1/\log C_9}}{\log C_9}. \]
Then if $|d| \geq C_g$, 
\[ 2^g \geq h(d) > \left(\frac{0.655}{\pi}\right)\frac{d^{1/2-1/\log C_9}}{\log C_9} 
\geq \left(\frac{0.655}{\pi}\right)\frac{C_g^{1/2-1/\log C_9}}{\log C_9} \geq 2^g, \]
a contradiction.  This completes the proof.
\end{proof}

We are now ready to prove main result of this section.

\begin{proof}[Proof of Proposition ${\ref{alldiscrim}}$]
We have already computed (in the previous section) that there are exactly $226$ such fundamental discriminants with $|d| \leq B$.  Therefore the proposition will follow from Lemma \ref{effconstant} and Lemma \ref{2224restrictsf} when it is shown that there are no fundamental discriminants $d<0$ with $\Cl(d)$ of type dividing $(2,\dots,2,4)$ satisfying one of the following conditions:
\begin{enumerate}
\item $4\cdot 67^4=B \leq |d| < C_9$; or
\item The integer $d$ has exactly $g$ distinct prime divisors, $g \in \{10,11\}$ and $C_9 \leq |d| < C_g$.
\end{enumerate}

Note that from the proof of Lemma \ref{effconstant}, we find $C_9=25593057435 \approx 2.5\cdot 10^{10}$, $C_{10}=116145031943 \approx 1.1 \cdot 10^{11}$, and $C_{11}=527083115400 \approx 5.2\cdot 10^{11}$.

The computations in (1) and (2) can be simplified by appealing to Lemma \ref{expc}: if $p \leq \sqrt[4]{|d|/4}$, then $(d/p)\neq 1$.  We then test for each prime $p$ such that $\sqrt[4]{|d|/4} < p \leq \sqrt{|d|/4}$ and $(d/p)=1$ if $\frakp^4$ is principal (working in the group of quadratic forms of discriminant $d$), where $(p)=\frakp\overline{\frakp}$.  To further rule out discriminants, we may also check given two such primes $p_1,p_2$ that $(\frakp_1\frakp_2)^2$ is principal.  For $d$ which satisfy all these conditions, we compute the class group $\Cl(d)$ itself (e.g.~using an algorithm of Shanks) and check explicitly if it is of type dividing $(2,\dots,2,4)$.  A computer search in \textsf{MAGMA} found no such $d$.  (The code is available from the author by request.)
\end{proof}

We also prove a complementary result which relies on a generalized Riemann hypothesis.

\begin{prop} \label{zetaD}
If the zeta function of the field $K=\bfQ(\sqrt{d})$ of discriminant $d<0$ does not have a zero in the interval $[1-(2/\log |d|),1)$ and the class group of $K$ is of type dividing $(2,\dots,2,4)$, then $|d|<1.3 \cdot 10^{10}$.
\end{prop}

\begin{prop}[Louboutin \cite{Louboutin}] \label{thmloub}
Let $K=\bfQ(\sqrt{d})$ be an imaginary quadratic field of discriminant $d$.  Suppose that the zeta function of $K$ does not have a zero in the interval $[1-(2/\log |d|),1)$.  Then 
\[ h(d) \geq \frac{\pi}{3e} \frac{\sqrt{|d|}}{\log |d|}, \]
where $e=\exp(1)$.
\end{prop}

\begin{proof}[Proof of Proposition $\ref{zetaD}$]
We follow \cite[Th\'eor\`eme 2]{Louboutin}.  Let $g$ be the number of distinct prime factors of the discriminant $d$.  Then $\#\Cl(d)[2]=2^{g-1}$ so $h(d) \leq 2^g$.  From Proposition \ref{thmloub}, we see that $2^g \geq (\pi/3e)\sqrt{|d|}/\log |d|$.  Recall that $|d| \geq \dhat_t=4p_1 \dots p_{t-1}$ whenever $d \neq -3$.  If we set
\[ t_0=\inf\left\{t \in \bfZ_{>0}: u \geq t \Rightarrow 2^u<\left(\frac{\pi}{3e}\right)\frac{(\dhat_u)^{1/2}}{\log \dhat_u}\right\}, \]
then $|d|<\dhat_{t_0}$ (see \cite{Louboutin}).  We compute easily that in this case $\dhat_{t_0}=4 \cdot 3 \cdot \ldots \cdot 29 < 1.3 \cdot 10^{10}$.  
\end{proof}

\begin{thm} \label{alldiscrimGRH}
Under the above Riemann hypothesis, there are exactly $226$ fundamental discriminants $d$ such that $\Cl(d)$ is of type dividing $(2,\dots,2,4)$, and $199$ such discriminants $D$ of nonmaximal orders.
\end{thm}

This follows from Proposition \ref{zetaD} and the computations performed in the proof of Proposition \ref{alldiscrim}.

\section*{Appendix: Ring Class Groups}

\setcounter{thm}{0}
\renewcommand{\thesection}{A}

In this appendix, we prove a proposition which characterizes ring class groups; we give a full statement for completeness.  

\begin{prop} \label{CRT}
For $f \in \bfZ_{>0}$, we have
\[ \frac{(A/fA)^*}{(\bfZ/f\bfZ)^*} \cong \prod_{p^e \parallel f} \frac{(A/p^e A)^*}{(\bfZ/p^e \bfZ)^*} \]
where $p$ is prime and $e>0$.  We have
\[ \frac{(A/2^e A)^*}{(\bfZ/2^e \bfZ)^*} \cong 
\begin{cases}
0, &\text{ if }d \equiv 1 \psmod{8}\text{ and }e=1; \\
\bfZ/2\bfZ \oplus \bfZ/2^{e-2} \bfZ, & \text{ if }d \equiv 1 \psmod{8}\text{ and }e \geq 2; \\
\bfZ/3\bfZ, &\text{ if }d \equiv 5 \psmod{8}\text{ and }e=1; \\
\bfZ/3\bfZ\oplus \bfZ/2\bfZ \oplus \bfZ/2^{e-2}\bfZ, &\text{ if }d \equiv 5 \psmod{8}\text{ and }e \geq 2; \\
\bfZ/2\bfZ \oplus \bfZ/2^{e-1}\bfZ, &\text{ if }d \equiv 4 \psmod{8}; \\
\bfZ/2^{e}\bfZ, &\text{ if }d \equiv 0\psmod{8};
\end{cases} \]
and
\[ \frac{(A/3^e A)^*}{(\bfZ/3^e \bfZ)^*} \cong 
\begin{cases}
\bfZ/2\bfZ \oplus \bfZ/3^{e-1}\bfZ, &\text{ if }d \equiv 1 \psmod{3}; \\
\bfZ/4\bfZ \oplus \bfZ/3^{e-1}\bfZ, &\text{ if }d \equiv 2 \psmod{3}; \\
\bfZ/3^{e}\bfZ, &\text{ if }d \equiv 3 \psmod{9}; \\
\bfZ/3\bfZ \oplus \bfZ/3^{e-1}\bfZ, &\text{ if }d \equiv 6 \psmod{9};
\end{cases} \]
and finally for $p \neq 2,3$, we have
\[ \frac{(A/p^e A)^*}{(\bfZ/p^e \bfZ)^*} \cong 
\begin{cases}
\bfZ/(p-1)\bfZ \oplus \bfZ/p^{e-1} \bfZ, &\text{ if }(d/p)=1; \\
\bfZ/(p+1)\bfZ \oplus \bfZ/p^{e-1} \bfZ, &\text{ if }(d/p)=-1; \\
\bfZ/p^e \bfZ, &\text{ if }(d/p)=0. \\
\end{cases} \]
\end{prop}

\begin{proof}
The first statement follows from the Chinese remainder theorem.  From Lemma \ref{exactseq}, we have
\[ \#\frac{(A/p^e A)^*}{(\bfZ/p^e \bfZ)^*}=p^e\left(1-\legen{d}{p}\frac{1}{p}\right). \]

We first treat the trivial case $p^e=2$: then $(\bfZ/2\bfZ)^*$ is the trivial group, and 
\[ (A/2A)^* \cong
\begin{cases}
\bfZ/3\bfZ, &\text{ if }(d/2)=-1; \\
\bfZ/2\bfZ, &\text{ if }(d/2)=0; \\
0, &\text{ if }(d/2)=1.
\end{cases} \]

Note that
\[ \frac{(A/p^e A)^*}{(\bfZ/p^e \bfZ)^*} \cong \frac{(A_p/p^e A_p)^*}{(\bfZ_p/p^e \bfZ_p)^*}, \]
where $A_p$ denotes the completion of $A$ at $p$ and $\bfZ_p$ the ring of $p$-adic integers.  
So if $(d/p)=1$, by \cite[\S II.5]{Neukirch} we have
\[ \frac{(A_p/p^eA_p)^*}{(\bfZ_p/p^e \bfZ_p)^*} \cong (\bfZ_p/p^e\bfZ_p)^* =
\begin{cases}
0, &\text{ if }p^e=2; \\
\bfZ/2\bfZ \oplus \bfZ/2^{e-2}\bfZ, &\text{ if $p=2$ and $e \geq 2$}; \\
\bfZ/(p-1)\bfZ \oplus \bfZ/p^{e-1}\bfZ, &\text{ otherwise}.
\end{cases} \]
From now on we assume $p^e \neq 2$ and $(d/p) \neq 1$.

Let $K_p$ denote the completion of $K$ at $p$, so that $A_p$ is its valuation ring with maximal ideal $\frakp$ and uniformizer $\pi$.  We denote by $v$ the unique valuation on $K_p$ normalized so that $v(p)=1$.
Let 
\[ V(A_p) = \{x\in A_p: v(x)>1/(p-1)\}. \]
It follows from \cite[Proposition II.5.4]{Neukirch} that there exists a (continuous) homomorphism $\log_p:A_p^*\to A_p$, with the property that $\log_p$ restricts to an isomorphism $1+V(A_p)\xrightarrow{\sim} V(A_p)$.  

One has an exact sequence
\[ 0 \to \frac{1+V(A_p)}{1+p^e A_p} \to \left(\frac{A_p}{p^e A_p}\right)^* \to \frac{A_p^*}{1+V(A_p)} \to 0. \]
We have an analogous exact sequence for $\bfZ_p$, and since $(1+V(A_p)) \cap \bfZ_p=1+V(\bfZ_p)$, it injects term-by-term into the one for $A_p$, yielding the following exact sequence:
\[ 
\tag*{$(\diamond)$}
0 \to \frac{\displaystyle{\frac{1+V(A_p)}{1+p^e A_p}}}{\displaystyle{\frac{1+V(\bfZ_p)}{1+p^e\bfZ_p}}} \to 
\frac{\displaystyle{\left(\frac{A_p}{p^eA_p}\right)^*}}{\displaystyle{\left(\frac{\bfZ_p}{p^e\bfZ_p}\right)^*}} \xrightarrow{\phi} \frac{\displaystyle{\frac{A_p^*}{1+V(A_p)}}}{\displaystyle{\frac{\bfZ_p^*}{1+V(\bfZ_p)}}} \to 0. \]

From the above, we see that by the logarithm map, 
\[ \frac{1+V(A_p)}{1+p^e A_p} \cong \frac{V(A_p)}{p^e A_p}\quad\text{ and }\quad 
\frac{1+V(\bfZ_p)}{1+p^e \bfZ_p} \cong \frac{V(\bfZ_p)}{p^e \bfZ_p}. \]  

Now let us assume that $p \neq 2,3$.  Then $V(A_p)=\frakp$, and $V(\bfZ_p)=p\bfZ_p$.  
We first analyze the group 
\[ \ker \phi=\frac{V(A_p)/p^e A_p}{V(\bfZ_p)/p^e \bfZ_p} \] 
from $(\diamond)$; we claim it is cyclic.  If $(d/p)=-1$, with $\eps \in A_p$ such that $A_p=\bfZ_p+\eps \bfZ_p$ as additive groups, then since $\frakp=V(A_p)=p\bfZ_p+p\eps\bfZ_p$, the element $p\eps$ generates the group $\ker \phi$.  If $(d/p)=0$, then $V(A_p)=\pi \bfZ_p+p \bfZ_p$, so $\pi$ is a generator.  It follows that
\[
\ker \phi \cong 
\begin{cases}
\bfZ/p^{e-1}\bfZ, &\text{ if $(d/p)=-1$}; \\
\bfZ/p^e\bfZ, &\text{ if $(d/p)=0$}.
\end{cases} \]

Now we analyze the image of $\phi$.  We have $A_p^*/(1+V(A_p)) \cong \mu(A_p)$ (see \cite[Proposition II.5.3]{Neukirch}), and this group can be computed as follows.  Since $[K_p:\bfQ_p]=2$ and the extension $\bfQ_p(\zeta_p)$ is a totally ramified extension of $\bfQ_p$ of degree $p-1$, we conclude that $A_p$ contains no $p$-power roots of unity.  Therefore
\[ \mu(A_p) \cong 
\begin{cases}
\bfZ/(p^2-1)\bfZ, &\text{ if $(d/p)=-1$}; \\
\bfZ/(p-1)\bfZ, &\text{ if $(d/p)=0$}.
\end{cases} \]
Since $\mu(\bfZ_p) \cong \bfZ/(p-1)\bfZ$, putting these two pieces together, we see that in the exact sequence $(\diamond)$, the kernel and image groups have orders which are relatively prime to each other and hence the exact sequence splits, and we obtain the result of the proposition.  

To conclude, we must treat the cases $p=2$, $p=3$.  Every field extension of $\bfQ_2$ of degree $2$ is isomorphic to $\bfQ_2(\sqrt{c})$ for $c \in \{-1,\pm 2,\pm 3,\pm 6\}$, and similarly for $\bfQ_3$ we have $c \in \{-1,-3,3,\}$.  We leave to the reader to verify the following: for $p=2$,
\[ \frac{(A/2^e A)^*}{(\bfZ/2^e\bfZ)^*} \cong
\begin{cases}
\la (-1+\sqrt{-3})/2 \ra \times \la \sqrt{-3} \ra \times \la 1+2\sqrt{-3} \ra, &\text{ if $c=-3$}; \\
\la \sqrt{-1} \ra \times \la 1+2\sqrt{-1} \ra, &\text{ if $c=-1$}; \\
\la 1+2\sqrt{c} \ra, &\text{ if $c=3$}; \\
\la 1+\sqrt{c} \ra, &\text{ if $2 \mid c$}; \\
\end{cases} \]
and for $p=3$,
\[ \frac{(A/3^e A)^*}{(\bfZ/3^e \bfZ)^*} \cong
\begin{cases}
\la 1+\sqrt{-1} \ra \times \la 1+3\sqrt{-1} \ra, &\text{ if $c=-1$}; \\
\la 1+\sqrt{3} \ra, &\text{ if $c=3$}; \\
\la (1+\sqrt{-3})/2 \ra \times \la 1+3\sqrt{-3} \ra, &\text{ if $c=-3$}.
\end{cases} \]
Computing the orders of these elements yields the conclusion of the proposition.
\end{proof}

\section*{Tables}

In Tables 1--2, we list equivalence classes with two fundamental discriminants ($\delta(C)=2$), then in Tables 3--5 those with three fundamental discriminants, then in Table 6 the exceptional cases with one fundamental discriminant (see Proposition \ref{k1eqk2}(iii)--(iv)).  Within each table, the classes are sorted by the smallest fundamental discriminant $d$ in each class.  Every form in an equivalence class has associated to it the same genus class field $P$ (Lemma \ref{samegenus}), denoted $\bfQ[a_1,\dots,a_r]=\bfQ(\sqrt{a_1},\dots,\sqrt{a_r})$.  The class group $\Cl_f(d)$ for each form is given by its type.  The set $E$ denotes the exceptional set for each equivalence class (\ref{except}).  

If $r \in \bfZ_{\geq 0}$, an abelian group $G$ is said to be \emph{of type $(\underbrace{2,\dots,2}_r,4)$} if 
\[ G \cong (\bfZ/2\bfZ)^r \oplus \bfZ/4\bfZ. \]
The group $G$ is \emph{of type dividing $(2,\dots,2,4)$} if there is an injection of groups
\[ G \hookrightarrow (\bfZ/2\bfZ)^r \oplus \bfZ/4\bfZ \]
for some $r \in \bfZ_{\geq 0}$.

In Tables 7--16, we list the orders of imaginary quadratic fields with class group of type dividing $(2,\dots,2,4)$, with at most possible exception (as in Theorem \ref{alldiscrim}).  In particular, there is no order with class group of type $(2,2,2,2,2)$ (unless this is the one exception!).  The tables are sorted by the isomorphism class of the class group, and within each table the classes are sorted by fundamental discriminant and then discriminant.

\begin{table}
\[ 
\begin{array}{cccc|cc|c}
Q & |D| & |d| & f & P & \Cl_f(d) & E \\
\hline \hline
\la 1, 0, 5 \ra & 20 & 20 & 1 & \bfQ[-1, 5] & (2) & \{5\} \\ 
\la 1, 0, 25 \ra & 100 & 4 & 5 & & (2) & \emptyset \\ 
\hline
\la 1, 0, 8 \ra & 32 & 8 & 2 & \bfQ[-1, 2] & (2) & \emptyset \\ 
\la 1, 0, 16 \ra & 64 & 4 & 4 & & (2) & \emptyset \\ 
\hline
\la 1, 0, 9 \ra & 36 & 4 & 3 & \bfQ[-1, -3] & (2) & \emptyset \\ 
\la 1, 0, 12 \ra & 48 & 3 & 4 & & (2) & \emptyset \\ 
\hline
\la 5, 0, 6 \ra & 120 & 120 & 1 & \bfQ[2, -3, 5] & (2, 2) & \{5\} \\ 
\la 11, 4, 14 \ra & 600 & 24 & 5 & & (2, 4) & \emptyset \\ 
\hline
\la 5, 0, 8 \ra & 160 & 40 & 2 & \bfQ[-1, 2, 5] & (2, 2) & \{5\} \\ 
\la 13, 8, 32 \ra & 1600 & 4 & 20 & & (2, 4) & \emptyset \\ 
\hline
\la 1, 0, 45 \ra & 180 & 20 & 3 & \bfQ[-1, -3, 5] & (2, 2) & \emptyset \\ 
\la 1, 0, 60 \ra & 240 & 15 & 4 & & (2, 2) & \emptyset \\ 
\hline
\la 5, 0, 9 \ra & 180 & 20 & 3 & \bfQ[-1, -3, 5] & (2, 2) & \{5\} \\ 
\la 9, 6, 26 \ra & 900 & 4 & 15 & & (2, 4) & \emptyset \\ 
\hline
\la 8, 0, 9 \ra & 288 & 8 & 6 & \bfQ[-1, 2, -3] & (2, 2) & \emptyset \\ 
\la 9, 6, 17 \ra & 576 & 4 & 12 & & (2, 4) & \emptyset \\ 
\hline
\la 1, 0, 120 \ra & 480 & 120 & 2 & \bfQ[-1, 2, -3, 5] & (2, 2, 2) & \emptyset \\ 
\la 1, 0, 240 \ra & 960 & 15 & 8 & & (2, 2, 2) & \emptyset \\ 
\hline
\la 5, 0, 24 \ra & 480 & 120 & 2 & \bfQ[-1, 2, -3, 5] & (2, 2, 2) & \{5\} \\ 
\la 21, 6, 29 \ra & 2400 & 24 & 10 & & (2, 2, 4) & \emptyset \\ 
\hline
\la 3, 0, 40 \ra & 480 & 120 & 2 & \bfQ[-1, 2, -3, 5] & (2, 2, 2) & \{3\} \\ 
\la 27, 12, 28 \ra & 2880 & 20 & 12 & & (2, 2, 4) & \emptyset \\ 
\hline
\la 3, 0, 56 \ra & 672 & 168 & 2 & \bfQ[-1, 2, -3, -7] & (2, 2, 2) & \{3\} \\ 
\la 20, 12, 27 \ra & 2016 & 56 & 6 & & (2, 2, 4) & \emptyset \\ 
\hline
\la 8, 0, 21 \ra & 672 & 168 & 2 & \bfQ[-1, 2, -3, -7] & (2, 2, 2) & \emptyset \\ 
\la 29, 12, 36 \ra & 4032 & 7 & 24 & & (2, 2, 4) & \emptyset \\ 
\hline
\la 3, 0, 80 \ra & 960 & 15 & 8 & \bfQ[-1, 2, -3, 5] & (2, 2, 2) & \{3\} \\ 
\la 27, 24, 32 \ra & 2880 & 20 & 12 & & (2, 2, 4) & \emptyset \\ 
\hline
\la 7, 6, 39 \ra & 1056 & 264 & 2 & \bfQ[-1, 2, -3, -11] & (2, 2, 4) & \emptyset \\ 
\la 7, 4, 76 \ra & 2112 & 132 & 4 & & (2, 2, 4) & \emptyset \\ 
\hline
\la 15, 12, 20 \ra & 1056 & 264 & 2 & \bfQ[-1, 2, -3, -11] & (2, 2, 4) & \emptyset \\ 
\la 23, 12, 36 \ra & 3168 & 88 & 6 & & (2, 2, 4) & \emptyset \\ 
\hline
\la 13, 6, 21 \ra & 1056 & 264 & 2 & \bfQ[-1, 2, -3, -11] & (2, 2, 4) & \emptyset \\ 
\la 13, 2, 61 \ra & 3168 & 88 & 6 & & (2, 2, 4) & \emptyset \\ 
\hline
\la 8, 0, 39 \ra & 1248 & 312 & 2 & \bfQ[-1, 2, -3, 13] & (2, 2, 2) & \emptyset \\ 
\la 15, 12, 44 \ra & 2496 & 39 & 8 & & (2, 2, 4) & \emptyset \\ 
\hline
\la 5, 4, 68 \ra & 1344 & 84 & 4 & \bfQ[-1, 2, -3, -7] & (2, 2, 4) & \emptyset \\ 
\la 5, 2, 101 \ra & 2016 & 56 & 6 & & (2, 2, 4) & \emptyset \\ 
\hline
\la 11, 8, 32 \ra & 1344 & 84 & 4 & \bfQ[-1, 2, -3, -7] & (2, 2, 4) & \emptyset \\ 
\la 11, 4, 92 \ra & 4032 & 7 & 24 & & (2, 2, 4) & \emptyset \\ 
\end{array} 
\]
\caption{Equivalence Classes $C$ of Forms ($\#\delta(C)=2$, $\#C=2$), 1~of~2} 
\end{table}

\setcounter{table}{0}

\begin{table}
\[ 
\begin{array}{cccc|cc|c}
Q & |D| & |d| & f & P & \Cl_f(d) & E \\
\hline \hline
\la 20, 4, 23 \ra & 1824 & 456 & 2 & \bfQ[-1, 2, -3, -19] & (2, 2, 4) & \emptyset \\ 
\la 23, 20, 44 \ra & 3648 & 228 & 4 & & (2, 2, 4) & \emptyset \\ 
\hline
\la 19, 4, 28 \ra & 2112 & 132 & 4 & \bfQ[-1, 2, -3, -11] & (2, 2, 4) & \emptyset \\ 
\la 19, 10, 43 \ra & 3168 & 88 & 6 & & (2, 2, 4) & \emptyset \\ 
\hline
\la 8, 0, 105 \ra & 3360 & 840 & 2 & \bfQ[-1, 2, -3, 5, -7] & (2, 2, 2, 2) & \emptyset \\ 
\la 32, 24, 57 \ra & 6720 & 420 & 4 & & (2, 2, 2, 4) & \emptyset \\ 
\hline
\la 21, 0, 40 \ra & 3360 & 840 & 2 & \bfQ[-1, 2, -3, 5, -7] & (2, 2, 2, 2) & \emptyset \\ 
\la 45, 30, 61 \ra & 10080 & 280 & 6 & & (2, 2, 2, 4) & \emptyset \\ 
\hline
\la 24, 0, 55 \ra & 5280 & 1320 & 2 & \bfQ[-1, 2, -3, 5, -11] & (2, 2, 2, 2) & \emptyset \\ 
\la 39, 36, 76 \ra & 10560 & 660 & 4 & & (2, 2, 2, 4) & \emptyset \\ 
\hline
\la 33, 0, 40 \ra & 5280 & 1320 & 2 & \bfQ[-1, 2, -3, 5, -11] & (2, 2, 2, 2) & \emptyset \\ 
\la 52, 36, 57 \ra & 10560 & 660 & 4 & & (2, 2, 2, 4) & \emptyset \\ 
\hline
\la 23, 4, 68 \ra & 6240 & 1560 & 2 & \bfQ[-1, 2, -3, 5, 13] & (2, 2, 2, 4) & \emptyset \\ 
\la 23, 18, 207 \ra & 18720 & 520 & 6 & & (2, 2, 2, 4) & \emptyset \\ 
\hline
\la 28, 12, 57 \ra & 6240 & 1560 & 2 & \bfQ[-1, 2, -3, 5, 13] & (2, 2, 2, 4) & \emptyset \\ 
\la 72, 48, 73 \ra & 18720 & 520 & 6 & & (2, 2, 2, 4) & \emptyset \\ 
\hline
\la 21, 12, 76 \ra & 6240 & 1560 & 2 & \bfQ[-1, 2, -3, 5, 13] & (2, 2, 2, 4) & \emptyset \\ 
\la 45, 30, 109 \ra & 18720 & 520 & 6 & & (2, 2, 2, 4) & \emptyset \\ 
\hline
\la 35, 30, 51 \ra & 6240 & 1560 & 2 & \bfQ[-1, 2, -3, 5, 13] & (2, 2, 2, 4) & \emptyset \\ 
\la 36, 12, 131 \ra & 18720 & 520 & 6 & & (2, 2, 2, 4) & \emptyset \\ 
\hline
\la 19, 14, 91 \ra & 6720 & 420 & 4 & \bfQ[-1, 2, -3, 5, -7] & (2, 2, 2, 4) & \emptyset \\ 
\la 19, 16, 136 \ra & 10080 & 280 & 6 & & (2, 2, 2, 4) & \emptyset \\ 
\hline
\la 28, 20, 85 \ra & 9120 & 2280 & 2 & \bfQ[-1, 2, -3, 5, -19] & (2, 2, 2, 4) & \emptyset \\ 
\la 45, 30, 157 \ra & 27360 & 760 & 6 & & (2, 2, 2, 4) & \emptyset \\ 
\hline
\la 51, 48, 56 \ra & 9120 & 2280 & 2 & \bfQ[-1, 2, -3, 5, -19] & (2, 2, 2, 4) & \emptyset \\ 
\la 59, 4, 116 \ra & 27360 & 760 & 6 & & (2, 2, 2, 4) & \emptyset \\ 
\hline
\la 33, 24, 88 \ra & 11040 & 2760 & 2 & \bfQ[-1, 2, -3, 5, -23] & (2, 2, 2, 4) & \emptyset \\ 
\la 57, 6, 97 \ra & 22080 & 1380 & 4 & & (2, 2, 2, 4) & \emptyset \\ 
\hline
\la 39, 6, 71 \ra & 11040 & 2760 & 2 & \bfQ[-1, 2, -3, 5, -23] & (2, 2, 2, 4) & \emptyset \\ 
\la 71, 70, 95 \ra & 22080 & 1380 & 4 & & (2, 2, 2, 4) & \emptyset \\ 
\hline
\la 76, 20, 145 \ra & 43680 & 10920 & 2 & \bfQ[-1, 2, -3, 5, -7, 13] & (2, 2, 2, 2, 4) & \emptyset \\ 
\la 96, 72, 241 \ra & 87360 & 5460 & 4 & & (2, 2, 2, 2, 4) & \emptyset \\ 
\hline
\la 88, 32, 127 \ra & 43680 & 10920 & 2 & \bfQ[-1, 2, -3, 5, -7, 13] & (2, 2, 2, 2, 4) & \emptyset \\ 
\la 127, 4, 172 \ra & 87360 & 5460 & 4 & & (2, 2, 2, 2, 4) & \emptyset \\ 
\hline
\la 57, 18, 193 \ra & 43680 & 10920 & 2 & \bfQ[-1, 2, -3, 5, -7, 13] & (2, 2, 2, 2, 4) & \emptyset \\ 
\la 148, 132, 177 \ra & 87360 & 5460 & 4 & & (2, 2, 2, 2, 4) & \emptyset \\ 
\hline
\la 55, 10, 199 \ra & 43680 & 10920 & 2 & \bfQ[-1, 2, -3, 5, -7, 13] & (2, 2, 2, 2, 4) & \emptyset \\ 
\la 159, 120, 160 \ra & 87360 & 5460 & 4 & & (2, 2, 2, 2, 4) & \emptyset \\ 
\end{array} 
\]
\caption{Equivalence Classes $C$ of Forms ($\#\delta(C)=2$, $\#C=2$), 2~of~2} 
\end{table}

\begin{table}
\[ 
\begin{array}{cccc|cc|c}
Q & |D| & |d| & f & P & \Cl_f(d) & E \\
\hline \hline
\la 1, 1, 4 \ra & 15 & 15 & 1 & \bfQ[-3, 5] & (2) & \emptyset \\ 
\la 1, 0, 15 \ra & 60 & 15 & 2 & & (2) & \emptyset \\ 
\la 1, 1, 19 \ra & 75 & 3 & 5 & & (2) & \emptyset \\ 
\hline
\la 2, 2, 11 \ra & 84 & 84 & 1 & \bfQ[-1, -3, -7] & (2, 2) & \{2\} \\ 
\la 8, 4, 11 \ra & 336 & 84 & 2 & & (2, 4) & \emptyset \\ 
\la 11, 2, 23 \ra & 1008 & 7 & 12 & & (2, 4) & \emptyset \\ 
\hline
\la 3, 0, 8 \ra & 96 & 24 & 2 & \bfQ[-1, 2, -3] & (2, 2) & \{3\} \\ 
\la 8, 8, 11 \ra & 288 & 8 & 6 & & (2, 2) & \emptyset \\ 
\la 11, 6, 27 \ra & 1152 & 8 & 12 & & (2, 4) & \emptyset \\ 
\hline
\la 5, 2, 5 \ra & 96 & 24 & 2 & \bfQ[-1, 2, -3] & (2, 2) & \emptyset \\ 
\la 5, 4, 20 \ra & 384 & 24 & 4 & & (2, 4) & \emptyset \\ 
\la 5, 2, 29 \ra & 576 & 4 & 12 & & (2, 4) & \emptyset \\ 
\hline
\la 7, 6, 7 \ra & 160 & 40 & 2 & \bfQ[-1, 2, 5] & (2, 2) & \emptyset \\ 
\la 7, 2, 23 \ra & 640 & 40 & 4 & & (2, 4) & \emptyset \\ 
\la 7, 4, 12 \ra & 320 & 20 & 4 & & (2, 4) & \emptyset \\ 
\hline
\la 2, 2, 23 \ra & 180 & 20 & 3 & \bfQ[-1, -3, 5] & (2, 2) & \{2\} \\ 
\la 8, 4, 23 \ra & 720 & 20 & 6 & & (2, 4) & \emptyset \\ 
\la 3, 0, 20 \ra & 240 & 15 & 4 & & (2, 2) & \{3\} \\ 
\hline
\la 3, 0, 16 \ra & 192 & 3 & 8 & \bfQ[-1, 2, -3] & (2, 2) & \{3\} \\ 
\la 4, 4, 19 \ra & 288 & 8 & 6 & & (2, 2) & \emptyset \\ 
\la 16, 8, 19 \ra & 1152 & 8 & 12 & & (2, 4) & \emptyset \\ 
\hline
\la 6, 6, 19 \ra & 420 & 420 & 1 & \bfQ[-1, -3, 5, -7] & (2, 2, 2) & \emptyset \\ 
\la 19, 12, 24 \ra & 1680 & 420 & 2 & & (2, 2, 4) & \emptyset \\ 
\la 19, 16, 31 \ra & 2100 & 84 & 5 & & (2, 2, 4) & \emptyset \\ 
\hline
\la 11, 8, 11 \ra & 420 & 420 & 1 & \bfQ[-1, -3, 5, -7] & (2, 2, 2) & \emptyset \\ 
\la 11, 6, 39 \ra & 1680 & 420 & 2 & & (2, 2, 4) & \emptyset \\ 
\la 11, 10, 50 \ra & 2100 & 84 & 5 & & (2, 2, 4) & \emptyset \\ 
\hline
\la 4, 4, 31 \ra & 480 & 120 & 2 & \bfQ[-1, 2, -3, 5] & (2, 2, 2) & \emptyset \\ 
\la 16, 8, 31 \ra & 1920 & 120 & 4 & & (2, 2, 4) & \emptyset \\ 
\la 15, 0, 16 \ra & 960 & 15 & 8 & & (2, 2, 2) & \emptyset \\ 
\hline
\la 12, 12, 13 \ra & 480 & 120 & 2 & \bfQ[-1, 2, -3, 5] & (2, 2, 2) & \emptyset \\ 
\la 13, 2, 37 \ra & 1920 & 120 & 4 & & (2, 2, 4) & \emptyset \\ 
\la 13, 4, 28 \ra & 1440 & 40 & 6 & & (2, 2, 4) & \emptyset \\ 
\hline
\la 12, 12, 17 \ra & 672 & 168 & 2 & \bfQ[-1, 2, -3, -7] & (2, 2, 2) & \emptyset \\ 
\la 17, 10, 41 \ra & 2688 & 168 & 4 & & (2, 2, 4) & \emptyset \\ 
\la 17, 4, 20 \ra & 1344 & 84 & 4 & & (2, 2, 4) & \emptyset \\ 
\hline
\end{array} 
\]
\caption{Equivalence Classes $C$ of Forms ($\#\delta(C)=2$, $\#C=3$), 1~of~2} 
\end{table}

\setcounter{table}{1}

\begin{table}
\[ 
\begin{array}{cccc|cc|c}
Q & |D| & |d| & f & P & \Cl_f(d) & E \\
\hline \hline
\la 13, 2, 13 \ra & 672 & 168 & 2 & \bfQ[-1, 2, -3, -7] & (2, 2, 2) & \emptyset \\ 
\la 13, 4, 52 \ra & 2688 & 168 & 4 & & (2, 2, 4) & \emptyset \\ 
\la 13, 8, 40 \ra & 2016 & 56 & 6 & & (2, 2, 4) & \emptyset \\ 
\hline
\la 8, 8, 41 \ra & 1248 & 312 & 2 & \bfQ[-1, 2, -3, 13] & (2, 2, 2) & \emptyset \\ 
\la 32, 16, 41 \ra & 4992 & 312 & 4 & & (2, 2, 4) & \emptyset \\ 
\la 20, 12, 33 \ra & 2496 & 39 & 8 & & (2, 2, 4) & \emptyset \\ 
\hline
\la 12, 12, 73 \ra & 3360 & 840 & 2 & \bfQ[-1, 2, -3, 5, -7] & (2, 2, 2, 2) & \emptyset \\ 
\la 48, 24, 73 \ra & 13440 & 840 & 4 & & (2, 2, 2, 4) & \emptyset \\ 
\la 33, 12, 52 \ra & 6720 & 420 & 4 & & (2, 2, 2, 4) & \emptyset \\ 
\hline
\la 31, 22, 31 \ra & 3360 & 840 & 2 & \bfQ[-1, 2, -3, 5, -7] & (2, 2, 2, 2) & \emptyset \\ 
\la 31, 18, 111 \ra & 13440 & 840 & 4 & & (2, 2, 2, 4) & \emptyset \\ 
\la 31, 10, 55 \ra & 6720 & 420 & 4 & & (2, 2, 2, 4) & \emptyset \\ 
\hline
\la 20, 20, 47 \ra & 3360 & 840 & 2 & \bfQ[-1, 2, -3, 5, -7] & (2, 2, 2, 2) & \emptyset \\ 
\la 47, 40, 80 \ra & 13440 & 840 & 4 & & (2, 2, 2, 4) & \emptyset \\ 
\la 47, 42, 63 \ra & 10080 & 280 & 6 & & (2, 2, 2, 4) & \emptyset \\ 
\hline
\la 28, 28, 37 \ra & 3360 & 840 & 2 & \bfQ[-1, 2, -3, 5, -7] & (2, 2, 2, 2) & \emptyset \\ 
\la 37, 18, 93 \ra & 13440 & 840 & 4 & & (2, 2, 2, 4) & \emptyset \\ 
\la 37, 24, 72 \ra & 10080 & 280 & 6 & & (2, 2, 2, 4) & \emptyset \\ 
\hline
\la 8, 8, 167 \ra & 5280 & 1320 & 2 & \bfQ[-1, 2, -3, 5, -11] & (2, 2, 2, 2) & \emptyset \\ 
\la 32, 16, 167 \ra & 21120 & 1320 & 4 & & (2, 2, 2, 4) & \emptyset \\ 
\la 32, 24, 87 \ra & 10560 & 660 & 4 & & (2, 2, 2, 4) & \emptyset \\ 
\hline
\la 41, 38, 41 \ra & 5280 & 1320 & 2 & \bfQ[-1, 2, -3, 5, -11] & (2, 2, 2, 2) & \emptyset \\ 
\la 41, 6, 129 \ra & 21120 & 1320 & 4 & & (2, 2, 2, 4) & \emptyset \\ 
\la 41, 10, 65 \ra & 10560 & 660 & 4 & & (2, 2, 2, 4) & \emptyset \\ 
\end{array} 
\]
\caption{Equivalence Classes $C$ of Forms ($\#\delta(C)=2$, $\#C=3$), 2~of~2} 
\end{table}

\begin{table}
\[ 
\begin{array}{cccc|cc|c}
Q & |D| & |d| & f & P & \Cl_f(d) & E \\
\hline \hline
\la 4, 4, 7 \ra & 96 & 24 & 2 & \bfQ[-1, 2, -3] & (2, 2) & \emptyset \\ 
\la 7, 6, 15 \ra & 384 & 24 & 4 & & (2, 4) & \emptyset \\ 
\la 7, 2, 7 \ra & 192 & 3 & 8 & & (2, 2) & \emptyset \\ 
\la 7, 4, 28 \ra & 768 & 3 & 16 & & (2, 4) & \emptyset \\ 
\hline
\la 8, 8, 17 \ra & 480 & 120 & 2 & \bfQ[-1, 2, -3, 5] & (2, 2, 2) & \emptyset \\ 
\la 17, 16, 32 \ra & 1920 & 120 & 4 & & (2, 2, 4) & \emptyset \\ 
\la 17, 14, 17 \ra & 960 & 15 & 8 & & (2, 2, 2) & \emptyset \\ 
\la 17, 6, 57 \ra & 3840 & 15 & 16 & & (2, 2, 4) & \emptyset \\ 
\end{array} 
\]
\caption{Equivalence Classes $C$ of Forms ($\#\delta(C)=2$, $\#C=4$)} 
\end{table}

\begin{table}
\[ 
\begin{array}{cccc|cc|c}
Q & |D| & |d| & f & P & \Cl_f(d) & E \\
\hline \hline
\la 1, 0, 24 \ra & 96 & 24 & 2 & \bfQ[-1, 2, -3] & (2, 2) & \emptyset \\ 
\la 1, 0, 48 \ra & 192 & 3 & 8 & & (2, 2) & \emptyset \\ 
\la 1, 0, 72 \ra & 288 & 8 & 6 & & (2, 2) & \emptyset \\ 
\hline
\la 7, 4, 52 \ra & 1440 & 40 & 6 & \bfQ[-1, 2, -3, 5] & (2, 2, 4) & \emptyset \\ 
\la 7, 6, 87 \ra & 2400 & 24 & 10 & & (2, 2, 4) & \emptyset \\ 
\la 7, 2, 103 \ra & 2880 & 20 & 12 & & (2, 2, 4) & \emptyset \\ 
\hline
\la 15, 0, 56 \ra & 3360 & 840 & 2 & \bfQ[-1, 2, -3, 5, -7] & (2, 2, 2, 2) & \emptyset \\ 
\la 39, 12, 44 \ra & 6720 & 420 & 4 & & (2, 2, 2, 4) & \emptyset \\ 
\la 36, 12, 71 \ra & 10080 & 280 & 6 & & (2, 2, 2, 4) & \emptyset \\ 
\end{array} 
\]
\caption{Equivalence Classes $C$ of Quadratic Forms ($\#\delta(C)=3$, $\#C=3$)} 
\end{table}

\begin{table}
\[ 
\begin{array}{cccc|cc|c}
Q & |D| & |d| & f & P & \Cl_f(d) & E \\
\hline \hline
\la 8, 0, 15 \ra & 480 & 120 & 2 & \bfQ[-1, 2, -3, 5] & (2, 2, 2) & \emptyset \\ 
\la 12, 12, 23 \ra & 960 & 15 & 8 & & (2, 2, 2) & \emptyset \\ 
\la 23, 22, 47 \ra & 3840 & 15 & 16 & & (2, 2, 4) & \emptyset \\ 
\la 23, 8, 32 \ra & 2880 & 20 & 12 & & (2, 2, 4) & \emptyset \\ 
\hline
\la 11, 2, 11 \ra & 480 & 120 & 2 & \bfQ[-1, 2, -3, 5] & (2, 2, 2) & \emptyset \\ 
\la 11, 4, 44 \ra & 1920 & 120 & 4 & & (2, 2, 4) & \emptyset \\ 
\la 11, 10, 35 \ra & 1440 & 40 & 6 & & (2, 2, 4) & \emptyset \\ 
\la 11, 8, 56 \ra & 2400 & 24 & 10 & & (2, 2, 4) & \emptyset \\ 
\hline
\la 8, 8, 23 \ra & 672 & 168 & 2 & \bfQ[-1, 2, -3, -7] & (2, 2, 2) & \emptyset \\ 
\la 23, 16, 32 \ra & 2688 & 168 & 4 & & (2, 2, 4) & \emptyset \\ 
\la 15, 6, 23 \ra & 1344 & 84 & 4 & & (2, 2, 4) & \emptyset \\ 
\la 23, 4, 44 \ra & 4032 & 7 & 24 & & (2, 2, 4) & \emptyset \\ 
\end{array} 
\]
\caption{Equivalence Classes $C$ of Quadratic Forms ($\#\delta(C)=3$, $\#C=4$)} 
\end{table}

\begin{table}
\[
\begin{array}{cccc|cc|c}
Q & |D| & |d| & f & P & \Cl_f(d) & E \\
\hline \hline
\la 1,1,1 \ra & 3 & 3 & 1 & \bfQ[-3] & (1) & \{3\} \\
\la 1,0,3 \ra & 12 & 3 & 2 & & (1) & \{3\} \\
\la 1,1,7 \ra & 27 & 3 & 3 & & (1) & \emptyset \\ 
\hline
\la 1,0,1 \ra & 4 & 4 & 1 & \bfQ[-1] & (1) & \{2\} \\
\la 1,0,4 \ra & 16 & 4 & 2 & & (1) & \emptyset \\
\end{array} 
\]
\caption{Equivalence Classes $C$ of Quadratic Forms ($\#\delta(C)=1$)} 
\end{table}

\begin{table}
\[
\begin{array}{ccc|ccc}
|d| & f & |D| & |d| & f & |D| \\
\hline \hline
3 & 1 & 3 & 8 & 1 & 8 \\
3 & 2 & 12 & 11 & 1 & 11 \\ 
3 & 3 & 27 & 19 & 1 & 19 \\ 
4 & 1 & 4 & 43 & 1 & 43 \\ 
4 & 2 & 16 & 67 & 1 & 67\\ 
7 & 1 & 7 & 163 & 1 & 163 \\ 
7 & 2 & 28 \\ 
\end{array}
\]
\caption{Orders of Quadratic Fields with Class Groups of Type $(1)$}
\end{table}

\begin{table}
\[
\begin{array}{ccc|ccc|ccc}
|d| & f & |D| & |d| & f & |D| & |d| & f & |D| \\ 
\hline \hline
3 & 4 & 48 & 15 & 1 & 15 & 115 & 1 & 115 \\
3 & 5 & 75 & 15 & 2 & 60 & 123 & 1 & 123 \\ 
3 & 7 & 147 & 20 & 1 & 20 & 148 & 1 & 148 \\ 
4 & 3 & 36 & 24 & 1 & 24 & 187 & 1 & 187 \\ 
4 & 4 & 64 & 35 & 1 & 35 & 232 & 1 & 232 \\ 
4 & 5 & 100 & 40 & 1 & 40 & 235 & 1 & 235 \\ 
7 & 4 & 112 & 51 & 1 & 51 & 267 & 1 & 267 \\ 
8 & 2 & 32 & 52 & 1 & 52 & 403 & 1 & 403 \\ 
8 & 3 & 72 & 88 & 1 & 88 & 427 & 1 & 427 \\ 
11 & 3 & 99 & 91 & 1 & 91 &  & &
\end{array}
\]
\caption{Orders of Quadratic Fields with Class Groups of Type $(2)$}
\end{table}

\begin{table}
\[
\begin{array}{ccc|ccc|ccc|ccc}
|d| & f & |D| & |d| & f & |D| & |d| & f & |D| & |d| & f & |D| \\ 
\hline \hline
3 & 11 & 363 & 39 & 1 & 39 & 184 & 1 & 184 & 723 & 1 & 723 \\
3 & 13 & 507 & 39 & 2 & 156 & 203 & 1 & 203 & 763 & 1 & 763 \\ 
4 & 6 & 144 & 43 & 3 & 387 & 219 & 1 & 219 & 772 & 1 & 772 \\
4 & 7 & 196 & 52 & 2 & 208 & 259 & 1 & 259 & 955 & 1 & 955 \\
4 & 8 & 256 & 55 & 1 & 55 & 291 & 1 & 291 & 1003 & 1 & 1003 \\
4 & 10 & 400 & 55 & 2 & 220 & 292 & 1 & 292 & 1027 & 1 & 1027 \\
7 & 3 & 63 & 56 & 1 & 56 & 323 & 1 & 323 & 1227 & 1 & 1227 \\
7 & 6 & 252 & 67 & 3 & 603 & 328 & 1 & 328 & 1243 & 1 & 1243 \\
8 & 4 & 128 & 68 & 1 & 68 & 355 & 1 & 355 & 1387 & 1 & 1387 \\
11 & 5 & 275 & 136 & 1 & 136 & 388 & 1 & 388 & 1411 & 1 & 1411 \\
19 & 3 & 171 & 148 & 2 & 592 & 568 & 1 & 568 & 1507 & 1 & 1507 \\
19 & 5 & 475 & 155 & 1 & 155 & 667 & 1 & 667 & 1555 & 1 & 1555 \\
20 & 2 & 80 & 163 & 3 & 1467
\end{array}
\]
\caption{Orders of Quadratic Fields with Class Groups of Type $(4)$}
\end{table}

\begin{table}
\[
\begin{array}{ccc|ccc|ccc}
|d| & f & |D| & |d| & f & |D| & |d| & f & |D| \\ 
\hline \hline
3 & 8 & 192 & 168 & 1 & 168 & 520 & 1 & 520 \\ 
7 & 8 & 448 & 195 & 1 & 195 & 532 & 1 & 532 \\ 
8 & 6 & 288 & 228 & 1 & 228 & 555 & 1 & 555 \\ 
15 & 4 & 240 & 232 & 2 & 928 & 595 & 1 & 595 \\ 
20 & 3 & 180 & 280 & 1 & 280 & 627 & 1 & 627 \\ 
24 & 2 & 96 & 312 & 1 & 312 & 708 & 1 & 708 \\ 
35 & 3 & 315 & 340 & 1 & 340 & 715 & 1 & 715 \\ 
40 & 2 & 160 & 372 & 1 & 372 & 760 & 1 & 760 \\ 
84 & 1 & 84 & 408 & 1 & 408 & 795 & 1 & 795 \\ 
88 & 2 & 352 & 435 & 1 & 435 & 1012 & 1 & 1012 \\ 
120 & 1 & 120 & 483 & 1 & 483 & 1435 & 1 & 1435 \\ 
132 & 1 & 132
\end{array}
\]
\caption{Orders of Quadratic Fields with Class Groups of Type $(2,2)$}
\end{table}

\begin{table}
\[
\begin{array}{ccc|ccc|ccc|ccc|ccc}
|d| & f & |D| & |d| & f & |D| & |d| & f & |D| & |d| & f & |D| & |d| & f & |D| \\
\hline \hline
3 & 16 & 768 & 84 & 2 & 336 & 308 & 1 & 308 & 987 & 1 & 987 & 2067 & 1 & 2067 \\
4 & 12 & 576 & 88 & 3 & 792 & 323 & 3 & 2907 & 1012 & 2 & 4048 & 2139 & 1 & 2139 \\ 
4 & 15 & 900 & 88 & 4 & 1408 & 328 & 2 & 1312 & 1032 & 1 & 1032 & 2163 & 1 & 2163 \\ 
4 & 20 & 1600 & 91 & 3 & 819 & 340 & 2 & 1360 & 1060 & 1 & 1060 & 2212 & 1 & 2212 \\ 
7 & 12 & 1008 & 91 & 5 & 2275 & 372 & 2 & 1488 & 1128 & 1 & 1128 & 2392 & 1 & 2392 \\ 
7 & 16 & 1792 & 115 & 3 & 1035 & 403 & 3 & 3627 & 1131 & 1 & 1131 & 2451 & 1 & 2451 \\ 
8 & 12 & 1152 & 132 & 2 & 528 & 427 & 3 & 3843 & 1204 & 1 & 1204 & 2632 & 1 & 2632 \\ 
11 & 15 & 2475 & 136 & 2 & 544 & 456 & 1 & 456 & 1240 & 1 & 1240 & 2667 & 1 & 2667 \\
20 & 4 & 320 & 148 & 3 & 1332 & 532 & 2 & 2128 & 1288 & 1 & 1288 & 2715 & 1 & 2715 \\
20 & 6 & 720 & 148 & 4 & 2368 & 552 & 1 & 552 & 1443 & 1 & 1443 & 2755 & 1 & 2755 \\
24 & 4 & 384 & 155 & 3 & 1395 & 564 & 1 & 564 & 1635 & 1 & 1635 & 2788 & 1 & 2788 \\
24 & 5 & 600 & 184 & 2 & 736 & 568 & 2 & 2272 & 1659 & 1 & 1659 & 2968 & 1 & 2968 \\
39 & 4 & 624 & 187 & 3 & 1683 & 580 & 1 & 580 & 1672 & 1 & 1672 & 3172 & 1 & 3172 \\
40 & 3 & 360 & 203 & 3 & 1827 & 616 & 1 & 616 & 1752 & 1 & 1752 & 3243 & 1 & 3243 \\
40 & 4 & 640 & 228 & 2 & 912 & 651 & 1 & 651 & 1768 & 1 & 1768 & 3355 & 1 & 3355 \\
51 & 5 & 1275 & 232 & 3 & 2088 & 708 & 2 & 2832 & 1771 & 1 & 1771 & 3507 & 1 & 3507 \\
52 & 3 & 468 & 232 & 4 & 3712 & 820 & 1 & 820 & 1780 & 1 & 1780 & 4123 & 1 & 4123 \\
52 & 4 & 832 & 235 & 3 & 2115 & 852 & 1 & 852 & 1947 & 1 & 1947 & 4323 & 1 & 4323 \\
55 & 4 & 880 & 260 & 1 & 260 & 868 & 1 & 868 & 1992 & 1 & 1992 & 5083 & 1 & 5083 \\
56 & 2 & 224 & 264 & 1 & 264 & 915 & 1 & 915 & 2020 & 1 & 2020 & 5467 & 1 & 5467 \\
56 & 3 & 504 & 276 & 1 & 276 & 952 & 1 & 952 & 2035 & 1 & 2035 & 6307 & 1 & 6307 \\
68 & 3 & 612
\end{array}
\]
\caption{Orders of Quadratic Fields with Class Groups of Type $(2,4)$}
\end{table}

\begin{table}
\[
\begin{array}{ccc|ccc}
|d| & f & |D| & |d| & f & |D| \\ 
\hline \hline
15 & 8 & 960 & 1092 & 1 & 1092 \\ 
120 & 2 & 480 & 1155 & 1 & 1155 \\ 
168 & 2 & 672 & 1320 & 1 & 1320 \\ 
280 & 2 & 1120 & 1380 & 1 & 1380 \\ 
312 & 2 & 1248 & 1428 & 1 & 1428 \\ 
408 & 2 & 1632 & 1540 & 1 & 1540 \\ 
420 & 1 & 420 & 1848 & 1 & 1848 \\ 
520 & 2 & 2080 & 1995 & 1 & 1995 \\ 
660 & 1 & 660 & 3003 & 1 & 3003 \\ 
760 & 2 & 3040 & 3315 & 1 & 3315 \\ 
840 & 1 & 840 \\ 
\end{array}
\]
\caption{Orders of Quadratic Fields with Class Groups of Type $(2,2,2)$}
\end{table}

\begin{table}
\[
\begin{array}{ccc|ccc|ccc|ccc|ccc}
|d| & f & |D| & |d| & f & |D| & |d| & f & |D| & |d| & f & |D| & |d| & f & |D| \\
\hline \hline
7 & 24 & 4032 & 372 & 4 & 5952 & 1288 & 2 & 5152 & 3432 & 1 & 3432 & 6708 & 1 & 6708 \\
15 & 16 & 3840 & 408 & 4 & 6528 & 1380 & 2 & 5520 & 3480 & 1 & 3480 & 6820 & 1 & 6820 \\ 
20 & 12 & 2880 & 420 & 2 & 1680 & 1428 & 2 & 5712 & 3588 & 1 & 3588 & 7035 & 1 & 7035 \\ 
24 & 10 & 2400 & 456 & 2 & 1824 & 1435 & 3 & 12915 & 3640 & 1 & 3640 & 7315 & 1 & 7315 \\ 
39 & 8 & 2496 & 520 & 3 & 4680 & 1540 & 2 & 6160 & 3795 & 1 & 3795 & 7395 & 1 & 7395 \\
40 & 6 & 1440 & 520 & 4 & 8320 & 1560 & 1 & 1560 & 3828 & 1 & 3828 & 7480 & 1 & 7480 \\
55 & 8 & 3520 & 532 & 3 & 4788 & 1672 & 2 & 6688 & 4020 & 1 & 4020 & 7540 & 1 & 7540 \\
56 & 6 & 2016 & 532 & 4 & 8512 & 1716 & 1 & 1716 & 4180 & 1 & 4180 & 7755 & 1 & 7755 \\
84 & 4 & 1344 & 552 & 2 & 2208 & 1752 & 2 & 7008 & 4260 & 1 & 4260 & 7995 & 1 & 7995 \\
84 & 5 & 2100 & 595 & 3 & 5355 & 1768 & 2 & 7072 & 4420 & 1 & 4420 & 8008 & 1 & 8008 \\
88 & 6 & 3168 & 616 & 2 & 2464 & 1860 & 1 & 1860 & 4440 & 1 & 4440 & 8052 & 1 & 8052 \\
120 & 4 & 1920 & 660 & 2 & 2640 & 1992 & 2 & 7968 & 4452 & 1 & 4452 & 8547 & 1 & 8547 \\
132 & 4 & 2112 & 708 & 4 & 11328 & 2040 & 1 & 2040 & 4488 & 1 & 4488 & 8680 & 1 & 8680 \\
168 & 4 & 2688 & 715 & 3 & 6435 & 2244 & 1 & 2244 & 4515 & 1 & 4515 & 8715 & 1 & 8715 \\
228 & 4 & 3648 & 760 & 3 & 6840 & 2280 & 1 & 2280 & 4740 & 1 & 4740 & 8835 & 1 & 8835 \\
232 & 6 & 8352 & 760 & 4 & 12160 & 2392 & 2 & 9568 & 5115 & 1 & 5115 & 8932 & 1 & 8932 \\
260 & 3 & 2340 & 952 & 2 & 3808 & 2436 & 1 & 2436 & 5160 & 1 & 5160 & 9867 & 1 & 9867 \\
264 & 2 & 1056 & 1012 & 3 & 9108 & 2580 & 1 & 2580 & 5187 & 1 & 5187 & 10948 & 1 & 10948 \\
280 & 3 & 2520 & 1012 & 4 & 16192 & 2632 & 2 & 10528 & 5208 & 1 & 5208 & 11067 & 1 & 11067 \\
280 & 4 & 4480 & 1032 & 2 & 4128 & 2760 & 1 & 2760 & 5412 & 1 & 5412 & 11715 & 1 & 11715 \\
308 & 3 & 2772 & 1092 & 2 & 4368 & 2968 & 2 & 11872 & 6195 & 1 & 6195 & 13195 & 1 & 13195 \\
312 & 4 & 4992 & 1128 & 2 & 4512 & 3108 & 1 & 3108 & 6420 & 1 & 6420 & 14763 & 1 & 14763 \\
340 & 3 & 3060 & 1140 & 1 & 1140 & 3192 & 1 & 3192 & 6580 & 1 & 6580 & 16555 & 1 & 16555 \\
340 & 4 & 5440 & 1240 & 2 & 4960 & 3220 & 1 & 3220 & 6612 & 1 & 6612 & & & \\
\end{array}
\]
\caption{Orders of Quadratic Fields with Class Groups of Type $(2,2,4)$}
\end{table}

\begin{table}
\[
\begin{array}{ccc}
|d| & f & |D| \\ 
\hline \hline
840 & 2 & 3360 \\ 
1320 & 2 & 5280 \\ 
1848 & 2 & 7392 \\ 
5460 & 1 & 5460 \\ 
\end{array}
\]
\caption{Orders of Quadratic Fields with Class Groups of Type $(2,2,2,2)$}
\end{table}

\begin{table}
\[
\begin{array}{ccc|ccc|ccc}
|d| & f & |D| & |d| & f & |D| & |d| & f & |D| \\
\hline \hline
280 & 6 & 10080 & 2280 & 2 & 9120 & 8680 & 2 & 34720 \\
420 & 4 & 6720 & 2760 & 2 & 11040 & 9240 & 1 & 9240 \\ 
520 & 6 & 18720 & 3192 & 2 & 12768 & 10920 & 1 & 10920 \\ 
660 & 4 & 10560 & 3432 & 2 & 13728 & 12180 & 1 & 12180 \\ 
760 & 6 & 27360 & 3480 & 2 & 13920 & 14280 & 1 & 14280 \\ 
840 & 4 & 13440 & 3640 & 2 & 14560 & 14820 & 1 & 14820 \\ 
1092 & 4 & 17472 & 4440 & 2 & 17760 & 17220 & 1 & 17220 \\ 
1320 & 4 & 21120 & 4488 & 2 & 17952 & 19320 & 1 & 19320 \\ 
1380 & 4 & 22080 & 5160 & 2 & 20640 & 19380 & 1 & 19380 \\ 
1428 & 4 & 22848 & 5208 & 2 & 20832 & 19635 & 1 & 19635 \\ 
1540 & 3 & 13860 & 5460 & 2 & 21840 & 20020 & 1 & 20020 \\
1540 & 4 & 24640 & 7140 & 1 & 7140 & 31395 & 1 & 31395 \\
1560 & 2 & 6240 & 7480 & 2 & 29920 & 33915 & 1 & 33915 \\
1848 & 4 & 29568 & 8008 & 2 & 32032 & 40755 & 1 & 40755 \\
2040 & 2 & 8160 & 8580 & 1 & 8580 \\ 
\end{array}
\]
\caption{Orders of Quadratic Fields with Class Groups of Type $(2,2,2,4)$}
\end{table}

\begin{table}
\[
\begin{array}{ccc}
|d| & f & |D| \\ 
\hline \hline
5460 & 4 & 87360 \\ 
9240 & 2 & 36960 \\ 
10920 & 2 & 43680 \\ 
14280 & 2 & 57120 \\ 
19320 & 2 & 77280 \\ 
\end{array}
\]
\caption{Orders of Quadratic Fields with Class Groups of Type $(2,2,2,2,4)$}
\end{table}

\end{document}